\documentclass[a4paper,12pt]{amsart}

\usepackage{hyperref}

\usepackage[top=3cm,left=2.5cm,right=2.5cm,bottom=3cm]{geometry}
\usepackage{relsize}
 \usepackage{lineno}

\usepackage{amsmath,amsthm,pb-diagram,amssymb,comment}
\usepackage{amsfonts,graphicx,color}
\usepackage{enumitem, fancyhdr, dsfont, multicol}
\usepackage{thmtools,cleveref}
\usepackage[normalem]{ulem}
\usepackage{stmaryrd}
\usepackage{textcomp}
\usepackage{mathtools}
\usepackage{fontawesome}
\usepackage{textcomp}

\usepackage{tikz,float}
\usepackage{mleftright}
\usepackage{mathrsfs}
\usepackage{tikz-cd}
\usepackage{dsfont}
\usepackage{stmaryrd}

\hypersetup{
    colorlinks=true,
    linkcolor=magenta,
    filecolor=magenta,
    urlcolor=magenta,
    citecolor=magenta,
}
\definecolor{caribbeangreen}{rgb}{0.0, 0.8, 0.6}

\setlist{topsep=0ex,itemsep=1ex}

    \DeclareMathOperator{\dom}{dom}

    \DeclareMathOperator{\ran}{ran}

    \newcommand{\Pwf}{\mathcal{P}}

    \DeclareMathOperator{\pts}{\mathcal{P}}

    \DeclareMathOperator{\cov}{cov}

    \newcommand{\Por}{\mathbb{P}}

    \newcommand{\Q}{\mathds{Q}}

    \DeclareMathOperator{\cf}{cf}

    \newcommand{\la}{\langle}
    \newcommand{\ra}{\rangle}

\DeclareMathOperator{\hgt}{\mathrm{ht}}

\newcommand{\Fn}{\mathrm{Fn}}

\newcommand{\cantor}{{}^\omega2}

\newcommand{\FAM}{\mathrm{FAM}}

\newcommand{\Cov}{\mathrm{Cov}}

\newcommand{\alt}{\mathrm{ht}}

\newcommand{\limit}{\mathrm{lim}}

\newcommand{\Int}{\mathrm{int}}

\newcommand{\p}{\mathrm{Pr}}

\newcommand{\E}{\mathrm{E}}
\newcommand{\Var}{\mathrm{Var}}

\newcommand{\Leb}{\mathrm{Leb}}

\newcommand{\Bin}{\mathrm{Bin}}

\newcommand{\Lev}{\mathrm{Lev}} 

\newcommand{\At}{\mathrm{At}}

\newcommand{\suc}{\mathrm{succ}}

\newcommand{\bfr}{\mathbf{r}}

\newcommand{\bbE}{\mathbb{E}}

\newcommand{\bbP}{\mathbb{P}}
\newcommand{\bbQ}{\mathbb{Q}}
\newcommand{\bbR}{\mathbb{R}}

\newcommand{\bbZ}{\mathbb{Z}}

\newcommand{\cB}{\mathscr{B}}

\newcommand{\gb}{\mathfrak{b}}

\newcommand{\calT}{\mathcal{T}}
\newcommand{\calC}{\mathcal{C}}

\newcommand{\calA}{\mathcal{A}}

\newcommand{\calB}{\mathcal{B}}
\newcommand{\calP}{\mathcal{P}}

\newcommand{\calN}{\mathcal{N}}
\newcommand{\calK}{\mathcal{K}}

\newcommand{\calL}{\mathcal{L}}

\newcommand{\calF}{\mathcal{F}}

\newcommand{\calY}{\mathcal{Y}}

\newcommand{\Ba}{\mathcal{B}\boldsymbol{a}} 

\newcommand{\finseq}{{}^{ < \omega} 2}
\newcommand{\finseqz}{{}^{ < \omega} Z}

\newcommand{\varp}{\varepsilon}

\newcommand{\rest}{{\restriction}}

\newenvironment{PROOF}[2][\proofname.]
   {\begin{proof}[#1]\renewcommand{\qedsymbol}{\textsquare$_{\mathrm{#2}}$}}
   {\end{proof}}

\newcommand{\concat}[2]{#1{}^{\frown}#2}
\newcommand{\power}[2]{ \! {}^{#1} \!#2}

\renewcommand{\setminus}{\smallsetminus}
\newcommand{\conj}{\mathrel{\mbox{\scriptsize $\wedge$}}}

\definecolor{sub0}{RGB}{29,32,137}
\definecolor{sub1}{RGB}{1,71,157}
\definecolor{sub2}{RGB}{1,104,183}
\definecolor{sub3}{RGB}{0,160,234}
\definecolor{sug}{RGB}{0,154,68}
\definecolor{suy}{RGB}{208,219,1}

\definecolor{redun}{rgb}{0.65, 0.11, 0.19}
\definecolor{greenun}{rgb}{0.58, 0.71, 0.23}
\definecolor{dodger}{rgb}{0.0,0.5,1.0}
\definecolor{carrotorange}{rgb}{0.93, 0.57, 0.13}

\thanks{This work was supported by the following grants: the Grant-in-Aid for Scientific Research (C)  23K03198, Japan Society for the Promotion of Science (first author), and the Austrian Science Fund (FWF): project number P33895 (second author).}

\author{Diego A.~Mej\'ia}
\address{Graduate School of System Informatics, Kobe University, 1-1 Rokkodai-cho, Nada-ku, Kobe, Hyogo 657-8501, Japan}
\email{damejiag@people.kobe-u.ac.jp}
\urladdr{\url{https://researchmap.jp/mejia?lang=en}}

\address{TU Wien, Faculty of Mathematics and Geoinformation, Institute of Discrete Mathematics and Geometry, Wiedner Hauptstrasse 8--10, A--1040 Vienna, Austria }
\email{\href{mailto:andres.zapata@tuwien.ac.at}{andres.zapata@tuwien.ac.at}}

\urladdr{\url{https://sites.google.com/view/andres-uribe-afuz/home}}


\date{\today}

\begin{document}

\makeatletter
\def\@roman#1{\romannumeral #1}
\makeatother

\newcounter{enuAlph}
\renewcommand{\theenuAlph}{\Alph{enuAlph}}

\theoremstyle{plain}
  \newtheorem{theorem}{Theorem}[section]
  \newtheorem{corollary}[theorem]{Corollary}
  \newtheorem{lemma}[theorem]{Lemma}
  \newtheorem{mainlemma}[theorem]{Main Lemma}
  \newtheorem{mainproblem}[theorem]{Main Problem}
  \newtheorem{construction}[theorem]{Construction}
  \newtheorem{prop}[theorem]{Proposition}
  \newtheorem{clm}[theorem]{Claim}
  \newtheorem{fact}[theorem]{Fact}
  \newtheorem{exer}[theorem]{Exercise}
  \newtheorem{question}[theorem]{Question}
  \newtheorem{problem}[theorem]{Problem}
  \newtheorem{cruciallem}[theorem]{Crucial Lemma}
  \newtheorem{conjecture}[theorem]{Conjecture}
  \newtheorem{assumption}[theorem]{Assumption}
  \newtheorem{teorema}[enuAlph]{Theorem}

  \newtheorem*{corolario}{Corollary}
\theoremstyle{definition}
  \newtheorem{definition}[theorem]{Definition}
  \newtheorem{example}[theorem]{Example}
  \newtheorem{remark}[theorem]{Remark}
    \newtheorem{hremark}[theorem]{Historical Remark}
    \newtheorem{observation}[theorem]{Observation}
    \newtheorem{notation}[theorem]{Notation}
  \newtheorem{context}[theorem]{Context}

\parindent 0pt

  \newtheorem*{defi}{Definition}
  \newtheorem*{acknowledgements}{Acknowledgements}

\numberwithin{equation}{theorem}
\renewcommand{\theequation}{\thetheorem.\arabic{equation}}

\def\sectionautorefname{Section}
\def\subsectionautorefname{Subsection}


\title{The measure algebra adding $\theta$-many random reals is $\theta$-FAM-linked}


\author{Andrés F. Uribe-Zapata}

\begin{abstract}
    The notion of $\theta$-FAM-linkedness, introduced in the second author's master thesis, is a formalization of the notion of strong FAM limits for intervals, whose initial form and applications have appeared in the work of Saharon Shelah, Jakob Kellner, and  Anda T\u{a}nasie, for controlling cardinals characteristics of the continuum in ccc forcing extensions. This generalization was successful in this thesis to establish
    a general theory of iterated forcing using finitely additive measures.
    
    In this paper, using probability theory tools developed in the same thesis, we refine Saharon Shelah's proof of the fact that random forcing is $\sigma$-FAM-linked and prove that any complete Boolean algebra with a strictly positive probability measure satisfying the $\theta$-density property is $\theta$-$\FAM$-linked. As a consequence, we get a new example of $\theta$-$\FAM$-linked forcing notions: the measure algebra adding $\theta$-many random reals. 
\end{abstract}

\maketitle


\section{Introduction}\label{SecIntro}

In 2000, Saharon Shelah in \cite{Sh00} built a forcing iteration using finitely additive measures to prove that, consistently, the covering number of the Lebesgue-null ideal, $\cov(\calN),$ may have countable cofinality. The iteration used random forcing, whose structure was fundamental both to be able to extend the iteration in each step and to apply it to the problem of the cofinality of $\cov(\calN).$ Later, in 2019, Jakob Kellner, Saharon Shelah, and Anda T\u anasie in \cite{KST} made new contributions to the method of iterations using finitely additive measures. In particular, they defined a new notion, called \emph{strong FAM limit for intervals} (see \cite[Def.~1.7]{KST}), which turned out to be the key to extending iterations using finitely additive measures at successor steps (see \cite[Lem.~2.25]{KST}, \cite[Thm.~4.3.16]{uribethesis}, and \cite[Sec.~7.2]{CMU}). Finally, in 2023, in the second author's master's thesis (see \cite[Def.~4.1.1]{uribethesis}), a new notion was introduced: the \emph{intersection number for forcing notions} (\autoref{s2}), which turned out to be the key to extending iterations using finitely additive measures at limit steps (see \cite[Thm.~4.3.18]{uribethesis} and \cite[Sec.~7.3]{CMU}).  Based on the ideas of the intersection number and strong FAM limits, in \cite{uribethesis} a new linkedness property called \emph{$\theta$-$\FAM$-linkedness} (\autoref{i10}), where $\theta$ is an infinite cardinal, was defined to present a general theory of iterated forcing using finitely additive measures (see \cite[Ch.~4]{uribethesis} and \cite[Sec.~5--7]{CMU}). In particular, such a property allows extending the iterations at limit and successor steps, preserving the unbounding number $\gb$, and preserving what in \cite{CMU} (a paper with Miguel A. Cardona where the iteration method is refined and further generalized) is defined as \emph{strongly-$\theta$-anti-Bendixson families}, which is essential to prove the consistency of $\cf(\cov(\calN)) = \aleph_{0}$. 

In \cite[Ex.~4.2.14]{uribethesis} (see also~\cite[Ex.~6.4]{CMU}), it was shown that every forcing notion $\bbP$ is $\vert \bbP \vert$-$\FAM$-linked, hence, for example, Cohen forcing is $\sigma$-$\FAM$-linked.\footnote{I.e.\ $\aleph_0$-$\FAM$-linked.} However, beyond this result, it is difficult to find more examples, particularly when $\theta$ is ``small''. In fact, for $\theta = \aleph_{0}$ only two more examples are known: random forcing and $\tilde{\bbE},$ a variant of the eventually different forcing notion $\bbE$ introduced in~\cite{KST} (see also~\cite{M24Anatomy}), which resembles a forcing notion by Haim Horowitz and Saharon Shelah (see \cite{HoroShe}). The proof that random forcing is $\sigma$-$\FAM$-linked appears implicitly in \cite[Lem.~2.17]{Sh00} and \cite[Lem.~2.18]{Sh00}, however, the proof of the first lemma is extremely difficult and, in the authors' opinion, has significant gaps. Concretely, it is not clear how the variance to which Saharon Shelah refers in the final part of the proof can be bounded. Additionally, the proof uses many probability theory tools that lack references, which are defined and developed in the second author's master's thesis (see \cite[Ch.~2]{uribethesis}) and in \cite{MU}.

When trying to clarify some points of the proof of \cite[Lem.~2.17]{Sh00} to obtain a complete and detailed proof that random forcing is $\sigma$-$\FAM$-linked, we noticed that it is possible to prove generalized versions of \cite[Lem.~2.17]{Sh00} and \cite[Lem.~2.18]{Sh00} (see \autoref{i19} and \autoref{i23}, respectively), which finally led us to a new and more general result that constitutes the main result of this paper:

\begin{teorema}\label{mt}
    Any complete Boolean algebra with a strictly positive probability measure satisfying the $\theta$-density property is $\theta$-$\FAM$-linked.
\end{teorema}

The density property of a finitely additive measure appears in~\autoref{m14}.

As a consequence, in \autoref{i99} we show a new example of $\theta$-$\FAM$-linked forcing notions for each $\theta$: the \emph{measure algebra adding $\theta$-many random reals} (see \autoref{2.3}). Consequently, in \autoref{i30}, we get a particular case due to Saharon Shelah: random forcing is $\sigma$-$\FAM$-linked (see \cite{Sh00} and \cite[Thm.~4.2.18]{uribethesis}).

\section{Preliminaries}

As usual, we denote by $\bbZ,$ $\bbQ$ and $\bbR$ the sets of integers, rational and real numbers, respectively, however, $\omega$ denotes the set of natural numbers. If $I \subseteq \bbR$ is an interval, we define $I_{\bbQ} \coloneqq I \cap \bbQ$.  Let $A, B$ be sets and $\alpha$ an ordinal number. We denote by ${}^{A} B$ the set of functions $f$ from $A$ into $B$ and ${}^{< \alpha} A \coloneqq \bigcup \{ {}^{\xi} A \colon \xi < \alpha  \}.$ Similarly, ${}^{\leq \alpha}A \coloneqq {}^{< \alpha +1 }A.$ For any $t \in {}^{< \alpha} A$ we define its \emph{length} by $\lg(t) \coloneqq \dom(t).$  We use the symbols ``$\langle$'' and ``$\rangle$'' to denote sequences and ``$\langle \ \rangle$'' to denote the empty sequence. If $A, B$ are non-empty, then $\Fn(A, B)$ is the set of \emph{finite partial functions} from $A$ into $B,$ that is, functions $f \colon X \to B,$ such that $X \subseteq A$ is finite.  If $B \subseteq A,$  we denote by $\chi_{B}$ the characteristic function of $B$ over $A.$ 

We assume the reader to be familiar with basic techniques of set theory (see e.g.~\cite{Kunen}), in particular with forcing theory. 



\subsection{Atoms in Boolean algebras} 

Let $\cB = (\cB, 0_{\cB}, 1_{\cB}, \wedge, \vee, \sim)$ be a Boolean algebra. We denote $\cB^{+} \coloneqq \cB \setminus \{ 0_{\cB} \}.$ For $B \subseteq \cB,$ the \emph{Boolean subalgebra generated by $B$}, denoted by $\langle B \rangle_{\cB}$, is the $\subseteq$-smallest Boolean subalgebra of $\cB$ containing $B.$ In this case, $B$ is called the \emph{generator set} of $\langle B \rangle_{\cB}.$ If $B$ is a finite set, we say that $\langle B \rangle_{\cB}$ is a \emph{finitely generated Boolean algebra}. It is well-known that $\cB$ can be endowed with a partial order: for any $a, b \in \cB,$ $ a \leq_{\cB} b $ iff $ a \wedge b = a.$ Also, we define $a \sim b \coloneqq a \, \wedge {\sim} b$. Notice that, for any $b \in \cB,$ $\cB_{\leq_{\cB} b} \coloneqq \{ a \in \cB \colon a \leq_{\cB} b \}$ is a Boolean algebra with $1_{\cB_{\leq_{\cB} b}} = b.$ Recall that, if $I$ is an ideal on $\cB,$ then we can define the quotient between $\cB$ and $I,$ which we denote by $\cB / I.$   An \emph{atom of $\cB$} is a $\leq_{\cB}$-minimal non-zero element of $\cB.$  Denote by $\At_{\cB}$ the set of all atoms of $\cB.$  We are particularly interested in the atomic structure of finitely generated Boolean algebras because we will be able to characterize their atoms in a particular way. For this, we introduce the following notation:

\begin{itemize}
    \item For any $b \in \cB$ and $d \in \{ 0, 1 \}$ we define:  
    $$b^{d} \coloneqq 
        \left\{ \begin{array}{lcc}
             b, &   \rm{if}  & d = 0, \\[1ex]
             {\sim}b, &  \rm{if}  & d = 1.
            \end{array}
        \right.$$ \index{$b^{d}$}

    \item For $B\subseteq \cB$ and $\sigma \in \Fn(B, 2),$ we set $$ a_{\sigma} \coloneqq \bigwedge_{b \in \dom(\sigma)} b^{\sigma(b)}. $$  
\end{itemize}

Notice that, for any $\sigma \in \Fn(B, 2), \, a_{\sigma} \in \cB.$

\begin{fact}\label{b50}
    Let $ \cB$ be a Boolean algebra generated by a finite set $B$. Then $$\At_{\cB} = \{ a_{\sigma} \colon \sigma \in {}^{B}2 \wedge a_{\sigma} \neq 0 \},$$ and it partitions $1_{\cB},$ that is, if $a, a' \in \At_{\cB}$ and $a \neq a',$ then $a \wedge a' = 0_{\cB},$ and  $1_{\cB} = \bigvee \At_{\cB}$. Furthermore, any $b \in \cB$ is partitioned by atoms.  
\end{fact}

\subsection{Finitely additive measures}\label{2.5} 

In this subsection, based on \cite[Ch.~3]{uribethesis} and  \cite{CMUP}, we present some basic notions about finitely additive measures. 


\begin{definition}
    Let $\cB$ be a Boolean algebra. A \emph{finitely additive measure on $\cB$} is a function $\Xi \colon \cB \to [0,\infty]$ satisfying the following conditions:   
    
\begin{enumerate}[label=\rm{(\arabic*)}]
    \item \label{m4a} $\Xi (0_{\cB})=0$,
            
    \item \label{m4b} $\Xi(a\vee b)=\Xi(a)+\Xi(b)$ whenever $a,b\in\cB$ and $a \wedge b= 0_{\cB}$.
\end{enumerate}

We say that $\Xi$ is a \emph{measure on $\cB$} if it satisfies~\ref{m4a} and
\begin{enumerate}
    \item[(2)'] $\Xi \! \left(\bigvee\limits_{n<\omega} b_n\right)=\sum\limits_{n<\omega} \Xi(b_n)$ whenever $\{ b_n \colon n < \omega \}  \subseteq \cB$, $\bigvee_{n<\omega}b_n$ exists, and $b_m\wedge b_n=0_\cB$ for $m\neq n$.
\end{enumerate}
        
    We exclude the trivial finitely additive measure, that is,  
    we will always assume $\Xi(1_\cB)>0$.  When $\Xi(b)>0$ for any $b \in \cB^{+}$ we say that $\Xi$ is \emph{strictly positive}.  Also, if $\Xi(1_\cB)=1,$ we say that $\Xi$ is a \emph{finitely additive probability measure}. 
\end{definition}

\begin{definition}
    Let $\cB$ be a Boolean algebra, $\Xi$ a finitely additive on it, and $b \in \cB$ with positive finite measure. We define the function 
    $\Xi_{b} \colon \cB \to [0, \infty]$ by $\Xi_{b}(a) \coloneqq \frac{\Xi(a \wedge b)}{\Xi(b)}$ for any $a \in \cB$.
\end{definition}

It is clear that $\Xi_{b}$ is a finitely additive probability measure. Also, when restricted to $\cB_{\leq_\cB b}$, it is a finitely additive probability measure on $\cB_{\leq_\cB b}$.

Recall that a \emph{field of sets over $K$} is a Boolean subalgebra of $\pts(K)$, where the latter is a Boolean algebra with the set-theoretic operations. In general, we will not be interested in finitely additive measures on fields of sets assigning a positive measure to finite sets. For this reason, we introduce the notion of \emph{free finitely additive measure}: 


\begin{definition}[{\cite[Def.~3.1.3]{uribethesis}}]\label{m12}
    If $K$ is a non-empty set, $\cB$ is a field of sets on $K$ and $\Xi$ is a finitely additive measure on $\cB$, we say that $\Xi$ is a \emph{free finitely additive measure} if, for any $k\in K$, $\{k\}\in\cB$ and $\Xi(\{k\}) = 0$. 
\end{definition}

It is clear that $\Xi_b$ is free whenever $\Xi$ is free and $b\in\cB$ has positive measure.

Based on \cite[Def.~5.4]{BCM2}, we introduce the notion of \emph{$\theta$-density property}: 


\begin{definition}\label{m14}
    Let $\cB$ a Boolean algebra, $\mu$ a strictly positive probability measure on $\cB,$ and $\theta$ an infinite cardinal. We say that $\mu$ satisfies the \emph{$\theta$-density property} if there exists an $S \subseteq \cB^{+}$ such that $\vert S \vert \leq \theta$ and, for any $\varp > 0$ and $b \in \cB^{+}$, there is some $s \in S$ such that $\mu_{s}(b) \geq 1 - \varp.$
\end{definition}

For example, any $\sigma$-centered Boolean algebra has a strictly positive finitely additive measure with the $\aleph_{0}$-density property (see \cite[Lem.~5.5]{BCM2}).

Now, we introduce integration over a field of sets. Fix a non-empty set $K$, a field of sets $\cB$ over $K$, and a finitely additive measure $\Xi \colon \cB \to [0, \infty).$ Motivated by the definition of Riemann's integral, if  $f \colon K \to \mathbb{R}$ is a bounded function we can naturally define $\int_{K} f d \Xi$ if it exists. In this case, we say that $f$ is \emph{$\Xi$-integrable}  (see \cite[Def.~3.5.3]{uribethesis} and \cite[Sec.~4.1]{CMUP}).  For example, any bounded function is $\Xi$-integrable when $\dom(\Xi) = \calP(K)$ (see \cite[Thm.~3.5.10]{uribethesis} and \cite[Sec.~4.2]{CMUP}). In fundamental aspects, the integral with respect to finitely additive measures behaves similarly to the Riemann integral, that is, we have available the basic properties of the integral such as those presented in \cite[Sec.~3.5]{uribethesis}. Next, we list some of these properties necessary for this work.

\begin{lemma}[{\cite[Thm.~3.5.12, Lem.~3.5.13 \& Lem.~3.5.14]{uribethesis}}]\label{t47}
    Let $f, g$ be $\Xi$-integrable functions and $c \in \bbR.$ Then: 

    \begin{enumerate}[label=\rm{(\arabic*)}]
        \item\label{t61} $fg$ and $cf$ are $\Xi$-integrable. In this case,  $\int_{K} (cf) d \Xi = c \int_{K} f d \Xi.$

        \item \label{t50} Let $\{ f_{i} \colon i < n \}$ a finite sequence of $\Xi$-integrable functions. Then $\sum_{i < n} f_{i}$ is $\Xi$-integrable and $$\int_{K} \left( \sum_{i < n} f_{i}\right) d \Xi = \sum_{i < n} \left( \int_{K} f_{i} \, d \Xi \right) \! .$$

        \item If $f \leq g$ then $ \int_{K} f d \Xi \leq \int_{K} g d \Xi.$
    \end{enumerate}
\end{lemma}

We also can integrate over subsets of $K$: for  a bounded function $f\colon K \to\bbR$, if $E \subseteq K$ and $\chi_{E} f$ is $\Xi$-integrable, we define $$\int_{E} f d \Xi \coloneqq \int_{K}  \chi_{E} f d \Xi.$$ 

\begin{lemma}[{\cite[Lem.~3.5.18 \& Thm.~3.5.22]{uribethesis}}]\label{t73}

    Let $E \subseteq K$ and $f\colon K \to\bbR$ be a bounded function. Then 
    
    \begin{enumerate}[label=\rm{(\arabic*)}]
        \item\label{t63} If $E \in \cB,$ then $\chi_{E}$ is $\Xi$-integrable and $\int_{K} \chi_{E}  d \Xi = \Xi(E).$ 

        \item\label{t63.0} If $n < \omega,$  $\langle E_{i} \colon i < n \rangle$ is a partition of $K$ into sets in $\dom(\Xi)=\cB$ and $\chi_{E_{i}} f$ is $\Xi$-integrable for any $i < n,$ then $f$ is $\Xi$-integrable and $$\int_{K} f d \Xi = \sum_{i < n} \left( \int_{E_{i}} f  d \Xi \right)\!.$$
    \end{enumerate}
\end{lemma}

If $E \in \cB$  and $f$ is $\Xi$-integrable, then by \autoref{t47}~\ref{t61}, 
$\chi_{E} f$ is $\Xi$-integrable. In general, this is the context in which we will use integration over subsets. 

The following result indicates that free finitely additive measures and their integrals can be approximated by the uniform measure on a finite set.

\begin{theorem}[{\cite[Thm.~3.5.26]{uribethesis} and \cite[Sec.~4.3]{CMUP}}]\label{t90} 
    Let $\Xi$ be a free finitely additive measure on $\calP(K)$ such that $\Xi(K) = \delta < \infty,$ $E \subseteq K$ with 
    $ \Xi(E) > 0$ and $I$ a non-empty finite set. For any $i \in I,$ let $f_{i} \colon K \to \bbR$ be a non-negative bounded function. Then, for all $\varp > 0$ and any $F \subseteq K,$ there exists a non-empty finite set $u \subseteq K \setminus F$ such that, for any $i \in I$, we have that: 
    $$ \left \vert \frac{\delta}{\vert u \vert}  \sum_{k \in u} f_{i}(k) -\frac{\delta}{\Xi(E)} \int_{E} f_{i} d \Xi    \right \vert <  \varp. $$ 
\end{theorem}

\subsection{Measure algebras}\label{2.3}

 In this subsection,  we review the \emph{measure algebra adding $\theta$-many random reals} for an infinite cardinal $\theta$, and we prove that the Lebesgue measure on it has the $\theta$-density property. This is a well-known fact, but we present its proof for completeness.

Recall that a \emph{measure algebra} is a pair $(\cB, \mu)$ where $\cB$ is a Boolean algebra and $\mu$ is a measure on $\cB.$ For example, random forcing $\calB(\cantor) / \calN(\cantor),$ where $\calB(\cantor)$ is the Borel $\sigma$-algebra on $\cantor$ and $\calN(\cantor)$ is the ideal of Lebesgue-null sets in $\cantor$,  with the Lebesgue measure, is a measure algebra. This type of algebras has been used as forcing notions to prove consistency results.

In a more general context than random forcing, Kenneth Kunen in~\cite{Ku84} presents an approach to deal with measure and category in the context of ${}^{I}2$ with $I$ infinite (hence possibly uncountable), endowed with the product topology of the discrete space $2=\{ 0, 1 \}$. Its topological basis is described by the basic clopen sets $[s] \coloneqq \{x \in {}^{I} 2 \colon  s \subseteq x \}$ for $s \in \Fn(I,2)$, and $\Ba({}^{I}2)$, the $\sigma$-algebra of \emph{Baire sets}, is defined as the $\sigma$-algebra generated by the clopen sets (which coincides with the product $\sigma$-algebra). We can also define the product measure $\Leb^{I}$ on the $\sigma$-algebra $\Ba({}^{I}2)$, where $\{ 0, 1 \}$ is endowed with the measure whose value on $\{ 0 \}$  and $\{ 1\}$  is $\frac{1}{2}$.  Then $\calN({}^{I} 2)$, the set of $\Leb^{I}$-null sets in ${}^{I} 2$, is an ideal on $\Ba({}^{I} 2)$ and, based on \cite[Def.~IV.7.33]{Kunen}, we can define: 


\begin{definition}\label{i13}
    For any infinite cardinal $\theta,$ the \emph{measure algebra adding $\theta$-many random reals} is the quotient $ \cB_{\theta} \coloneqq  \Ba({}^{\theta} 2) / \calN({}^{\theta} 2).$ 
\end{definition} 

Abusing the notation, for any infinite cardinal $\theta,$ we define $\Leb^{\theta} \colon \cB_{\theta} \to [0, \infty]$ by $\Leb^{\theta}([A]_{\calN({}^{\theta} 2)}) \coloneqq \Leb^{\theta}(A)$. Then:   

\begin{theorem}\label{i3.0}
    For any infinite cardinal $\theta,$ $\Leb^{\theta}$ on $\cB_{\theta}$  is a strictly positive probability measure satisfying the $\theta$-density property. 
\end{theorem}

\begin{PROOF}[\textbf{Proof}]{\autoref{i3.0}}
    It is clear that $\Leb^{\theta}$ is a strictly positive probability measure on $\cB_{\theta},$ hence we concentrate on the $\theta$-density property. Consider $S$ as the set of clopen sets in ${}^{\theta}2,$ and let $\varp > 0$. By compactness, we have that $\vert S \vert = \theta.$ Let $A \in \Ba({}^{\theta} 2)$ be such that $\Leb^{\theta}(A) >0.$ Then, we can find a countable set $I \subseteq \theta$ such that $A = B \times {}^{\theta \setminus I} 2$ for some $B \in \Ba({}^{I} 2).$ Notice that, by Fubini's theorem, $\Leb^{I}(B) = \Leb^{\theta}(A) > 0.$ On the other hand,  since $I$ is countable, by Lebesgue's density theorem (see \cite[Thm.~3.20]{Oxtoby}), there exists a clopen set $C^{\bullet}$ in ${}^{I} 2$ such that  $\Leb^{I}_{C^{\bullet}}(B) \geq 1 - \varp.$ Define $\pi_{I} \colon {}^{\theta} 2 \to {}^{I} 2$ such that, for any $x \in {}^{\theta} 2,$ $\pi_{I}(x) \coloneqq x \, {\rest} \, I,$ and consider $C \coloneqq \pi_{I}^{-1}[C^{\bullet}].$ Therefore $C \in S,$ moreover $C = C^{\bullet} \times {}^{\theta \setminus I} 2.$ Notice that, again by Fubini's theorem, $\Leb^{\theta}(A \cap C) = \Leb^{I}(B \cap C^{\bullet})$ and $\Leb^{\theta}(C) = \Leb^{I}(C^{\bullet}).$ As a consequence, we have that, $$ \Leb^{\theta}_{C}(A) = \frac{\Leb^{\theta}(A \cap C)}{\Leb^{\theta}(C)} = \frac{\Leb^{I}(B \cap C^{\bullet})}{\Leb^{I}(C^{\bullet})} = \Leb^{I}_{C^{\bullet}}(B) \geq 1 - \varp.$$ 

    Thus, $\Leb^{\theta}$ on $\cB_{\theta}$ has the $\theta$-density property. 
\end{PROOF}

\subsection{Basic probability theory}

The notions defined in this subsection are classical notions of basic probability theory as in, e.g.~\cite{ross98}. We say that $ \Omega \coloneqq (\Omega, \calA, \p)$ is a \emph{probability space} if $\Omega$ is a non-empty set, $\calA$ is a $\sigma$-algebra on $\Omega$ and $\p \colon \calA \to [0, 1]$ is a measure such that $\p(\Omega) = 1.$ In this case, we say that $\p$ is a \emph{probability measure on $\Omega$}. Elements in $\calA$ are called \emph{events} and, if $E, F \in \calA,$ then $\p(E)$ is called the \emph{probability of success of $E$}. 

In the context of probability theory, all the probability spaces that we are going to consider in this paper are finite and discrete. So the following lemma will be very useful to provide any finite set with a probability space structure:

\begin{fact}\label{p6}
    Let $\Omega$ be a finite set and $\p \colon \Omega \to [0, 1]$ such that $\sum_{o \in \Omega} \p(o) = 1.$ Then there exists a (unique) probability function $\p_{\Omega} \colon \calP(\Omega) \to [0, 1]$ such that $(\Omega, \calP(\Omega), \p_{\Omega})$ is a probability space and, for any $o \in \Omega, \, \p_{\Omega}(\{ o \}) = \p(o).$ 
\end{fact}



It will be useful not to make a distinction between $\p$ and $\p_{\Omega}$ in the Lemma above. That is, to ease the notation, if $o \in \Omega$ and $\p \colon \calP(\Omega) \to [0, 1]$ then, for any $o \in \Omega,$ we denote $\p(o) \coloneqq \p(\{ o \}).$ This justifies not making any distinction between the functions $\p$ and $\p_{\Omega}$ in \autoref{p6}, that is, when we want to define a probability space over a finite set $\Omega$, we are going to define a function $\p \colon \Omega \to [0, 1]$ that satisfies the conditions of \autoref{p6} and denote $\p_{\Omega}$ by $\p.$

Fix, for the rest of this subsection, a probability space $  (\Omega, \calA, \p).$ We recall the notion of \emph{random variable} on $\Omega$: 

\begin{definition}\label{p11}
    We say that a function $X \colon \Omega \to \bbR$ is a \emph{random variable on $\Omega$} if, for any $a \in \bbR, \, \{ o \in \Omega \colon X(o) \leq a \} \in \calA.$ 
    
\end{definition}

Notice that, in measure theoretic terms, a random variable is simply an $\calA$-measurable function. 


We say that a \emph{trial} is an experiment where there are only two possible outcomes, one with probability $p$ and the other with probability $1-p$. If $X$ counts the number of successes in a sequence of $n$ independent trials, each with probability of success $p,$ we say that $X$ has  \emph{binomial distribution with parameters $n, p$}, and we denote it as $X \sim \Bin(n, p).$ 


By the definition of a random variable, for any $x \in \bbR, \, \p(\{ o \in \Omega \colon X(o) = x \}), \,  \p(\{ o \in \Omega \colon X(o) \leq  x \})$ and $\p(\{ o \in \Omega \colon X(o) \geq x \})$  are defined. So, to simplify the writing, we denote $ \p[X = x] \coloneqq \p(\{ o \in \Omega \colon X(o) = x \}), $  $ \p[X \leq x] \coloneqq \p(\{ o \in \Omega \colon X(o) \leq  x \})$, and similarly for $ \p[X \geq x].$

Now, we recall the expected value, the variance, and the covariance for random variables: 

\begin{definition}\label{p30}
    Let $X, Y$ be random variables on $\Omega.$ Then:  

    \begin{enumerate}[label=\rm{(\arabic*)}]
        \item $\E[X] \coloneqq \sum_{r \in \ran(X)} r \cdot \p[X = r]$ is called the \emph{expected value} of $X.$ \index{expected value} \index{$\E[X]$}

        \item $\Cov[X, Y] \coloneqq \E[XY] - \E[X] \cdot \E[Y]$ is called the \emph{covariance} of $X$ and $Y.$ \index{covariance} 

        \item $\Var[X] \coloneqq \Cov[X, X]$ is called the \emph{variance} of $X.$ \index{variance}

    \end{enumerate}
\end{definition}

The expected value and the variance when $X$ has a binomial distribution is well-known: if $X \sim \Bin(n, p),$ then $\E[X] = np$ and $\Var[X] = np(1-p).$




        




Now, we review some elementary properties of the variance and covariance: 

\begin{theorem}\label{p39}
    Let $X$ and $Y$ be random variables on $\Omega$ such that $\vert \E[X] \vert < \infty,$ and $r \in \bbR.$ Then,  

    \begin{enumerate}[label=\rm{(\arabic*)}]
        \item\label{p39.1} $\Cov[X, Y]$ is a bilinear function. 

        \item\label{p39.2} $\Cov[X, r] = \Cov[r, Y] = 0.$
        
        \item\label{p39.3} If $\langle X_{n} \colon n < n^{\ast} \rangle$ is a sequence of random variables on $\Omega,$ then $$\Var \! \left[ \sum_{i < n } a_{i} \, X_{i} \right]   = \sum_{i < n}a_{i}^{2} \, \Var[X_{i}] + \sum_{i, j \leq n, \, i \neq j} a_{i} a_{j} \Cov[X_{i}, X_{j}].$$
    \end{enumerate}
\end{theorem}

Moreover, we have that $\Var[X] \geq 0$ and $\Var[X + r] = \Var[X].$ 

Finally, we state a result that will be fundamental in this work: the one known as \emph{Chebyshev's inequality}: 

\begin{theorem}\label{p42}
    Let $X$ be a random variable on $\Omega$ with $ \vert \E[X] \vert < \infty.$ Then for all  $\varp > 0,$ $\p[ \, \vert X - \E[X] \vert \geq \varp] \leq \frac{\Var[X]}{\varp^{2}}.$
\end{theorem}

\subsection{Trees}



Let $Z$ be a non-empty set. A \emph{subtree of ${}^{<\omega}Z$} is a set $\calT \subseteq \finseqz$ such that, for any $\rho, \eta \in \finseqz,$ if $ \rho \subseteq \eta$ and $\eta \in \calT,$ then $\rho \in \calT.$ For instance, for any $n^{\ast} < \omega$, ${}^{n^{\ast} \geq }2$ is a subtree of $\finseq$ called the \emph{complete binary tree} of height $n^{\ast} +1.$  

Fix a tree $\calT$ and $\rho, \eta \in \calT.$ Elements in $\calT$ are called \emph{nodes}, and if $\rho \subseteq \eta,$ we say that $\eta$ is an \emph{extension of $\rho$}. We say that $\rho$ and $\eta$ are \emph{compatible} if they are compatible as functions.  The \emph{height of $\rho$ in $\calT$} is $\alt_{\calT}(\rho) \coloneqq \dom(\rho)$ and  the \emph{height of $\calT$}  is $\alt(\calT) \coloneqq \sup \{ \alt_{\calT}(\rho) + 1 \colon \rho \in \calT \}.$ When the context is clear, we  simply denote $\alt_{\calT}(\rho)$ as $\alt(\rho)$. An \emph{infinite branch} of $\calT$ is an element of $z \in {}^{\omega}Z$ such that, for any $n < \omega,$ $z \, {\rest} \, n \in \calT.$  The set of infinite branches of $\calT$ is denoted by $[\calT].$ Also, for any $\rho \in \finseqz,$  we define $[\rho] \coloneqq \{x  \in {}^{\omega} Z \colon \rho \subseteq  x \}. $ 

Below we define other notions related to trees: 

\begin{definition}\label{y8}
    Let $\calT$ be a subtree of $\finseqz$ and $\rho \in \calT.$ We define:

    \begin{enumerate}[label=\rm{(\arabic*)}]
        \item For any $h < \omega,$ the \emph{$h$-th level of $\calT$} is  $\Lev_{h}(\calT) \coloneqq \calT \cap {}^{h} Z,$ 

        \item $\calT_{\geq \rho} \coloneqq \{ \eta \in \calT \colon \rho \subseteq \eta \}$ is the set of \emph{successors of $\rho$} in $\calT.$ 

        \item $\suc_{\calT}(\rho) \coloneqq  \calT_{\geq \rho} \cap \Lev_{\alt_{\calT}(\rho) + 1}(\calT),$ that is, $\suc_{\calT}(\rho)$ is the set of \emph{immediate successors of $\rho$} in $\calT.$ 

        \item $\max(\calT) \coloneqq \{ \rho \in \calT \colon \suc_{\calT}(\rho) = \emptyset \},$ that is, it is the set of \emph{maximal nodes} of $\calT.$ 


    \end{enumerate}

    When the context is clear we simply write  ``$\suc(\rho)$'' instead of ``$\suc_{\calT}(\rho)$''.
\end{definition}

Notice that we have, for any $\rho \in \calT,$ $\alt(\rho) < \omega$ and, if $[\calT] \neq \emptyset,$ then $\alt(\calT) = \omega.$  Also, $\rho \in \Lev_{j}(\calT)$ if, and only if, $\alt(\rho) = j.$
On the other hand, it is clear that $\calT_{\geq \rho}$ is isomorphic to some subtree of ${}^{< \omega} Z$  (see \autoref{f40}). Therefore, from now on, without loss of generality, we can consider it as a tree. 

\begin{definition}\label{y11}
    Let $\calT$ a subtree of $\finseqz.$  We say that $\calT$ is \emph{well-pruned tree} of height $n^{\ast} < \omega$ if $\calT \neq \emptyset$ and $\max(\calT) = \Lev_{n^{\ast}}(\calT).$  Similarly, $\calT$ is a \emph{well-pruned tree} of height $\omega$ if $\calT \neq \emptyset$ and, for any $\rho \in \calT, \, \suc_{\calT}(\rho) \neq \emptyset$. 
\end{definition}

    


\subsection{Probability trees}

In this subsection, based on \cite[Ch.~2]{uribethesis} and \cite{MU}, we present the necessary notions and results about \emph{probability trees}.  

\begin{definition}[{\cite[Def.~2.3.1]{uribethesis}}]\label{p45}
    Let $Z$ be a non-empty set. Say that a subtree $\calT$ of ${}^{<\omega}Z$ is a \emph{probability tree} if it is a well-pruned subtree of ${}^{<\omega} Z$ with associated probability spaces $\la \suc(\rho),\calA_{\rho}^{\calT},\p^\calT_\rho\ra$
    for any $\rho \in \calT \setminus \max(\calT)$, 
    where $[\suc(\rho)]^{< \aleph_{0}} \subseteq \calA_{\rho}^{\calT}$.
\end{definition} 

As we mentioned earlier, we can consider $\calT_{\geq \rho}$ as a tree. Furthermore, it inherits the probability structure from $\calT$ in a natural way:

\begin{lemma}\label{p47}
    Any probability tree $\calT$ induces a probability structure on $\calT_{\geq \rho}$ for any $\rho \in \calT$.   
\end{lemma}

\begin{PROOF}{\autoref{p47}}
    Notice that, if $\eta \in \calT_{\geq \rho} \setminus \max(\calT_{\geq \rho}),$ then $\eta \in \calT \setminus \max(\calT)$ and $\suc_{\calT}(\eta) = \suc_{\calT_{\geq \rho}}(\eta),$ so we can use the same probability $\sigma$-algebra and measure of $\suc_\calT(\eta)$ for $\suc_{\calT_{\geq\rho}}(\eta)$.
    %
\end{PROOF}

Since probability trees have a probability space structure at the immediate successors of each non-maximal node, a probability space structure is naturally induced at each level of the tree when levels are finite (even countable):

\begin{theorem}[{\cite[Thm.~2.3.2]{uribethesis} and \cite[Sec.~4.1]{MU}}]\label{p49}
    Every probability tree $\calT$ with finite levels induces a (discrete) probability space on each of its levels. 
\end{theorem}

For the rest of this subsection, we assume that each probability tree we consider has finite levels. 

If $\calT$ is a probability tree and $n<\hgt(\calT)$, we denote by $\Pr_{\Lev_n(\calT)}$ the probability measure determined at level $n$ of $\calT$. The following result allows us to decompose the probability of the successors of $\rho$ at the level $h+n$ of $\calT,$ in terms of the probability at the level $n$ of $\calT_{\geq \rho}$.

\begin{lemma}[{\cite[Def.~2.3.7]{uribethesis} and \cite[Sec.~4.4]{MU}}]\label{p62}
    Let $\calT$ be a probability subtree of $ \finseqz,$ $h < \omega$ and $\rho \in \Lev_{h}(\calT).$ Let $0 < n < \omega$ and $\eta \in \Lev_{n}(\calT_{\geq \rho}).$  Then, $$\p_{\Lev_{h+n}(\calT)}(\eta) = \p_{\Lev_{n}(\calT_{\geq \rho})}(\eta) \cdot \p_{\Lev_{h}(\calT)}(\rho).$$
\end{lemma}

The proof of \autoref{i19}, the key point to prove the main result in this paper, is quite long and technical. To simplify it a bit, we use the \emph{relative expected value} in probability trees as was introduced in \cite{uribethesis} and \cite{MU}: 



\begin{definition}[{\cite[Def.~2.3.5]{uribethesis} and \cite[Sec.~4.4]{MU}}]\label{p60}
     Let $\calT$ be a probability tree on $ \finseqz,$ $h < \omega$ and $\rho \in \Lev_{h}(\calT).$  Let $0 < n < \omega$ and $X$ be a random variable on $\Lev_{h + n}(\calT).$ Then, we define: $$\E_{\Lev_{h+n}(\calT)}[X \colon \eta\, \rest\, h = \rho] \coloneqq \E_{\Lev_{n}(\calT_{\geq \rho})}[X \,  {\rest} \, \Lev_{n}(\calT_{\geq \rho})],$$ and call it the \emph{relative expected value of $X$ with respect to $\rho$ in $\Lev_{h+n}(\calT)$}. Notice that the $\eta$ above is a dummy variable, i.e.\ the expected value of $X$ is calculated by varying $\eta$ over the nodes in $\calT$ at level $h+n$ that extend $\rho$. 

     When the context is clear, we simply write $E_{h+n}[X \colon \eta \, {\rest} \, h = \rho]$ or even $\E[X \colon \eta \rest h = \rho]$ instead of $\E_{\Lev_{h+n}(\calT)}[X \colon \eta \rest h = \rho].$  \index{$E_{h+n}[X \colon \eta \, {\rest} \, h = \rho]$} \index{$\E[X \colon \eta \rest h = \rho]$}
\end{definition}

Since the relative expected value is defined in terms of a usual expected value, it is clear that: 

\begin{theorem}[{\cite[Thm.~2.3.6]{uribethesis}}]\label{p61}
      Let $\calT$ be a probability tree on $ \finseqz, \, h < \omega$ and $\rho \in \Lev_{h}(\calT).$ Consider $0 < n < \omega,$ two random variables $X, Y$ on $\Lev_{h + n}(\calT)$ and $r, s \in \bbR.$ Then,  $\E_{h + n}[r X + s Y \colon \eta \rest h = \rho] = r \E_{h + n}[X \colon \eta \rest h = \rho] + s \E_{h+n}[Y \colon \eta \rest h = \rho].$
\end{theorem}

To calculate a relative expected value, we can use the \emph{intermediate steps}. 
\begin{theorem}[{\cite[Thm.~2.3.8]{uribethesis} and \cite[Sec.~4.4]{MU}}]\label{p64}
    Let $\calT$ be a probability subtree of $ \finseqz$, $ \rho \in \Lev_{h}(\calT)$  where $h < \omega,$ and let  $0 < n < m < \omega$. If $X$ is a random variable on $\Lev_{h+m}(\calT),$ then  $\E_{h+m}[X \colon \nu \rest h = \rho] = \E_{h+n}[ \,  \E_{h+m}[ X \colon \nu \rest (h +n) = \eta ] \colon \eta \rest h = \rho ].$ (See \autoref{f40}.) \vspace{-0.4cm}
\end{theorem}
  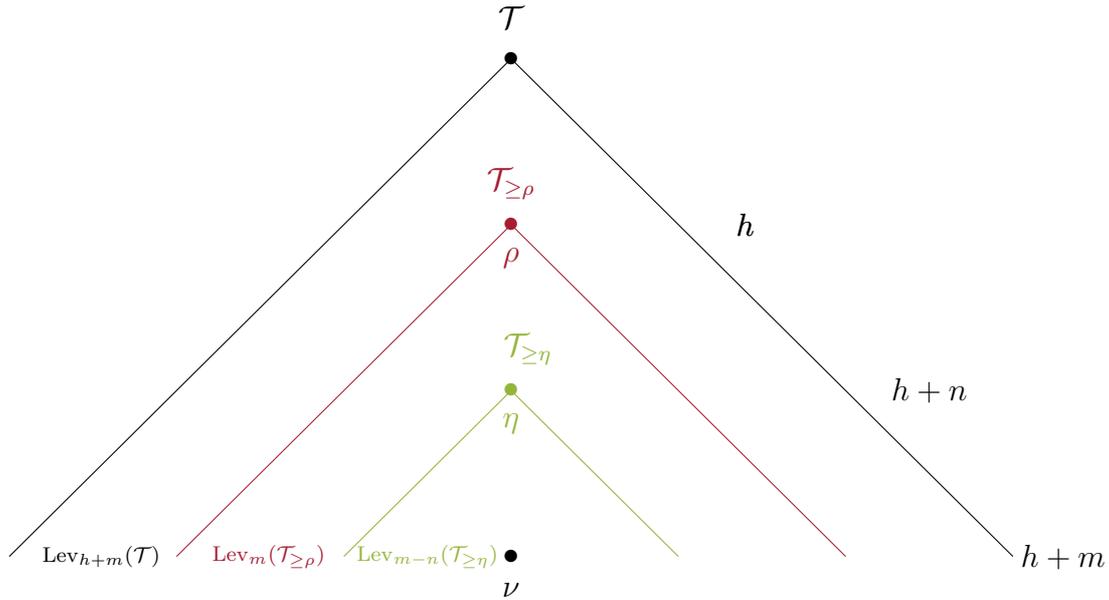
\begin{figure}[ht]
        \centering
        \begin{tikzpicture}[scale=1.1]

            \node at (6, 6) {$\bullet$}; 
            \node[redun] at (6, 4) {$\bullet$}; 
            \node[greenun] at (6, 2) {$\bullet$}; 
            \node at (6, 0) {$\bullet$}; 


            \draw (6, 6) -- (0, 0);
            \draw (6, 6) -- (12, 0);


            \draw[redun] (6, 4) -- (2, 0);
            \draw[redun] (6, 4) -- (10, 0);


            \draw[greenun] (6, 2) -- (4, 0);
            \draw[greenun] (6, 2) -- (8, 0);


            \node at (8.8, 4) {$h$};
            \node at (11, 2) {$h+ n$};
            \node at (12.6, 0) {$h + m$};


            \node at (6, 6.5) {$\calT$};
            \node[redun] at (6, 4.5){$\calT_{\geq \rho}$};
            \node[greenun] at (6.2, 2.5) {$\calT_{\geq \eta}$};
            \node at (8.8, 4) {$h$};         


            \node at (1.1, 0) {\tiny $\Lev_{h+m}(\calT)$};
            
            \node[redun] at (3.1, 0) {\tiny $\Lev_{m}(\calT_{\geq \rho})$};
            
            \node[greenun] at (5, 0) {\tiny $\Lev_{m-n}(\calT_{\geq \eta})$}; 


            \node at (6, -0.4) {$\nu$};
            \node[greenun] at (6, 1.6) {$\eta$};
            \node[redun] at (6, 3.6) {$\rho$};
        \end{tikzpicture}
        
        \caption{A graphic example of the situation in \autoref{p64}.}
        \label{f40}
    \end{figure}    

Finally, as a consequence,  we can express the expected value of $X$ in terms of the relative expected value:

\begin{corollary}[{\cite[Cor.~2.3.9]{uribethesis} and \cite[Sec.~4.4]{MU}}]\label{p66}
    Let $\calT$ be a probability subtree of $ \finseqz$ and $0 < n < m < \omega$. If $X$ is a random variable on $\Lev_{m}(\calT),$ then $\E_{m}[X] = \E_{n}[\, \E_{m}[X \colon \nu \rest n = \eta]].$ 
\end{corollary}

\subsection{$\mu$-$\FAM$-linkedness}

In this section, we introduce the \emph{intersection number for forcing notions} and the linkedness property called $(\Xi, \bar{P}, \varp)$-linkedness, the necessary definitions to define the notion of \emph{$\theta$-$\FAM$-linkedness}. For an in-depth and detailed study of such a linkedness property see \cite[Sec.~4.2 \& Sec.~4.3]{uribethesis} and \cite[Sec.~5]{CMU}; and for the notion of \emph{intersection number}, see \cite[Sec.~4.1]{uribethesis} and \cite{IntNumU24}. 


\begin{definition}[{\cite[Def.~4.1.1]{uribethesis}}]\label{s2}
    Let $\bbP$ be a forcing notion and $Q \subseteq \bbP.$ 

    \begin{enumerate}[label=\rm{(\arabic*)}]
        \item For a finite sequence $\bar{q} = \langle q_{i} \colon i < n \rangle \in {}^{n}\bbP,$ we define $$ i_{\ast}^{\bbP}(\bar{q}) \coloneqq \max \{ \vert F \vert \colon F \subseteq n \conj  \{ q_{i} \colon i \in F \} \ \text{has a lower bound in} \ \bbP  \}. $$ 

        \item  The \emph{intersection number of $Q$ in $\bbP$}, denoted by $\rm{int}^{\bbP}(Q),$ is defined by $$\Int^{\bbP}(Q) \coloneqq \inf \left \{ \frac{i_{\ast}^{\bbP}(\bar{q})}{n} \colon \bar{q} \in Q^{n}  \conj   n \in \omega \smallsetminus \{ 0 \}  \right\}.$$
    \end{enumerate}
\end{definition}

The intersection number was originally defined for Boolean algebras and has the following connection with measures.

\begin{lemma}[{\cite[Prop.~1]{Kelley59}}]\label{s9}
    Let $\cB$ be a Boolean algebra and $\Xi$ a finitely additive probability measure on $\cB$. Consider $\bbP \coloneqq \cB^{+}$ and $\delta \in [0, 1].$ If $Q \coloneqq \{ p \in \bbP \colon \Xi(p) \geq \delta \},$ then $\Int^{\bbP}(Q) \geq \delta.$
\end{lemma}


Although originally this result was due to John Kelley, a recent and complete proof can be found in \cite[Lem.~4.1.7]{uribethesis} and \cite{IntNumU24}. 

\begin{definition}[{\cite[Sec.~5.1]{CMU}}]\label{i1}
    Let $K$ and $W$ be non-empty sets, $\Xi$ a finitely additive probability measure on $\Pwf(K)$, 
    $\bar{P} = \langle \, P_{k} \colon k \in K \rangle$ a partition of $W$ into finite non-empty sets, $\varp_{0} \in (0, 1)$ and let $\bbP$ be a forcing notion. Say that $Q \subseteq \bbP$ is \emph{$(\Xi, \bar{P}, \varp_{0})$-linked} if there are a $\Por$-name of a finitely additive measure $\dot \Xi_Q$ on $\pts(K)$ extending $\Xi$ and a function $\limit^{\Xi} \colon \power{W}{Q} \to \bbP$ such that, for any $\bar q = \la q_\ell\colon \ell\in W\ra \in {}^W \! Q$, 
    \[\limit^{\Xi}\bar q\Vdash \text{``} \int_K \frac{|\{\ell \in P_k\colon q_\ell \in \dot G\}|}{|P_k|}d\dot \Xi_Q \geq 1-\varp_0 \text{''}.\]
\end{definition}

The following characterization of ``$(\Xi,\bar P,\varp_0)$-linked'' was first proved in \cite[Thm.~4.2.5]{uribethesis} when $\Xi$ is free and fixing $\Xi^-$ as the uniform measure on $|u|$, considering that this notion was firstly defined in~\cite[Def.~4.2.2]{uribethesis} in the form of~\ref{i2II}, which is closer to the original approach from \cite{Sh00,KST}.\footnote{In \cite[Def.~1.10]{KST}, this notion is called \emph{strong-$\FAM$-limit for intervals.}} However, this approach only makes sense when $\Xi$ is free (or with the \emph{uniform approximation property (uap)} as defined in \cite{CMUP}, see also \cite[Sec.~3]{CMU}). Later, with Cardona~\cite{CMU}, we obtained the equivalence below regardless of whether $\Xi$ is free (or has the uap). This is the reason why, contrary to tradition, we adopt \autoref{i1} as the main definition of ``$(\Xi,\bar P,\varp_0)$-linked''.

\begin{theorem}[{\cite[Sec.~5.5]{CMU}}]\label{i2}
    Fix $K,\ W,\ \Xi,\ \bar P,\ \varp_0,\ \bbP,$ and $Q$ as in \autoref{i1} and let\/ $\limit^{\Xi} \colon \power{W}{Q} \to \bbP$.\footnote{In \cite{KST}, this function is called \emph{$\FAM$-limit}.} Then, the following statements are equivalent.
    \begin{enumerate}[label = \rm (\Roman*)]
        \item\label{i2I} $Q$ is $(\Xi,\bar P,\varp_0)$-linked in $\Por$ witnessed by $\limit^\Xi$.
        \item\label{i2II} Given
        \begin{itemize}
        \item $i^{\ast} < \omega$ and $ \bar{q}^{i} = \langle q_{\ell}^{i} \colon \ell \in W \rangle \in \power{W}{Q},$ for each $i < i^{\ast},$

        \item $m^{\ast} < \omega$ and a partition  $\langle B_{m} \colon m < m^{\ast} \rangle$ of $K,$

        \item $\varp > 0$ and $q \in \bbP$ such that, for all $i < i^{\ast},$ $ q \leq \limit^{\Xi}(\bar{q}^{i}),$
        \end{itemize}
        there are a non-empty finite set $u\subseteq K$, a probability measure $\Xi^-$ on $\pts(u)$ and some $q'\leq q$ in $\bbP$ such that
    \begin{enumerate}[label=\rm{(\arabic*')}]
        \item\label{i2a} $\left| \Xi^-(u \cap B_{m}) - \Xi(B_{m}) \right| < \varp,$ for all $m < m^{\ast},$ and

        \item\label{i2b} $\sum_{k \in u} \frac{\vert \{ \ell \in P_{k} \colon  q' \leq q_{\ell}^{i} \} \vert}{\vert P_{k} \vert} \Xi^-(\{k\}) \geq 1 - \varp_{0} - \varp,$ for all $i < i^{\ast}.$
    \end{enumerate}
    \end{enumerate}
    Moreover, if $\Xi$ is free, $\Xi^-$ in~\ref{i2II} can be found as the uniform measure on $u$, i.e.\ each point in $u$ has measure $\frac{1}{\vert u \vert}$.
\end{theorem}

Finally, we can introduce the notion of $\theta$-$\FAM$-linkedness: 

\begin{definition}[{\cite[Def.~4.2.8]{uribethesis}} and~{\cite[Sec.~6]{CMU}}]\label{i10}
    Let $\bbP$ be a forcing notion, and $\theta$ be a cardinal. Denote by $\calY_*$ the class of all pairs $(\Xi,\bar P)$ such that $\Xi$ is a finitely additive probability measure on some $\pts(K)$ with $K$ non-empty, and $\bar P = \la P_k \colon k\in K\ra$ is a partition of some set $W$ into non-empty finite sets. Different pairs in $\calY_*$ may be associated with different $K$ and $W$, and note that $K= K_{\bar P}\coloneqq\dom \bar P$ and $W= W_{\bar P}\coloneqq\bigcup_{k\in K}P_k$.
    \begin{enumerate}[label=\rm{(\arabic*)}]
        \item\label{i10.1} Let $\calY\subseteq\calY_*$ be a class. 
        We say that $\bbP$ is \emph{$\theta$-$\calY$-linked}, if there exists a sequence $\langle Q_{\alpha, \varp} \colon \alpha < \theta \conj  \varp \in (0, 1)_{\bbQ} \rangle$ of subsets of $\bbP$, such that:
        
        \begin{enumerate}[label=\rm{(\alph*)}]
            \item\label{i10.1.a} for any $(\Xi,\bar P)\in\calY$, each  $Q_{\alpha, \varp}$ is $(\Xi, \bar{P}, \varp)$-linked,
    
            \item\label{i10.1.b} for every $\varp \in (0, 1)_{\bbQ}, \, \bigcup_{\alpha < \theta} Q_{\alpha, \varp}$ is dense in $\bbP,$
    
    
            \item\label{i10.1.c} for any $\alpha < \theta$ and $\varp \in (0, 1)_{\bbQ} \,, \Int(Q_{\alpha, \varp}) \geq 1 - \varp.$ 
        \end{enumerate}

        \item We say that $\bbP$ is \emph{$\theta$-$\FAM$-linked} if it is $\theta$-$\calY_\omega$-linked, where $\calY_\omega\coloneqq \{ (\Xi,\bar P)\in \calY_* \colon K_{\bar P} = \omega \text{ and }\Xi \text{ is free}\}$.
    \end{enumerate}

    When $\theta = \aleph_{0},$ we write ``$\sigma$-$\calY$-linked'' instead of ``$\aleph_{0}$-$\calY$-linked'', likewise for ``$\sigma$-$\FAM$-linked''.
\end{definition}

\begin{remark}\label{rm-i10}
    In \cite[Sec.~6]{CMU}, we proved that~\ref{i10.1.c} in \autoref{i10}~\ref{i10.1} is redundant when $\calY$ contains some $(\Xi',\bar P')$ satisfying that $\Xi'(\{k\in K_{\bar P'} \colon |P'_k|=n\}) = 0$ for all $n<\omega$ (which implies that $\Xi'$ is free). This holds, for example, when $K_{\bar P'} = \omega$, $\Xi'$ is free and $\lim_{k\to\infty} |P'_k| = \infty$.
\end{remark}

We even deal with a uniform version of \autoref{i10}.

\begin{definition}\label{i3}
    Let $\calY \subseteq \calY_*$ be a class, $\theta$ a cardinal, and $\Por$ a forcing notion. We say that $\Por$ is \emph{uniformly $\theta$-$\calY$-linked} if it is $\theta$-$\calY$-linked, witnessed by some sequence $\langle Q_{\alpha, \varp} \colon \alpha < \theta \conj  \varp \in (0, 1)_{\bbQ} \rangle$ of subsets of $\Por$, satisfying the following additional requirement for any $(\Xi,\bar P)\in \calY$: there is a $\Por$-name $\dot \Xi^*$ of a finitely additive measure extending $\Xi$ and, for each $\alpha<\theta$ and $\varp\in(0,1)_\bbQ$, there is a function $\limit^{\alpha,\varp}\colon {}^{W_{\bar P}} Q_{\alpha,\varp}\to \Por$  such that, for any $\bar q\in {}^{W_{\bar P}} Q_{\alpha,\varp}$, 
    \[\limit^{\alpha,\varp}\bar q \Vdash \text{``} \int_{K_{\bar P}} \frac{|\{\ell \in P_k\colon q_\ell \in \dot G\}|}{|P_k|}d\dot \Xi^* \geq 1-\varp \text{''}.\]
\end{definition}

The difference between \autoref{i10} and \autoref{i3} is that, in the former, the name $\dot\Xi^*$ of the finitely additive measure extending $\Xi$ witnessing the linkedness property of any $Q_{\alpha,\varp}$ depends on this set, while in \autoref{i3}, $\dot \Xi^*$ does not depend on $Q_{\alpha,\varp}$ (but the limit function still does).

We have a characterization of \autoref{i3} similar to \autoref{i2}.

\begin{theorem}[{\cite[Thm.~6.16]{CMU}}, cf.~{\cite[Thm.~3.10]{M24Anatomy}}]\label{i4}
    Let $\calY\subseteq\calY^*$ be a class, $\Por$ a forcing notion and $\langle Q_{\alpha, \varp} \colon \alpha < \theta \conj  \varp \in (0, 1)_{\bbQ} \rangle$ a sequence of subsets of $\Por$ satisfying~\ref{i10.1.b} and~\ref{i10.1.c} of \autoref{i10}. Then, the following statements are equivalent. 
    \begin{enumerate}[label = \rm (\Roman*)]
        \item\label{i4I} $\Por$ is uniformly $\theta$-$\calY$-linked witnessed by $\langle Q_{\alpha, \varp} \colon \alpha < \theta \conj  \varp \in (0, 1)_{\bbQ} \rangle$.

        \item\label{i4II} For each $(\Xi,\bar P)\in \calY$, $\alpha<\theta$ and $\varp\in (0,1)_\bbQ$, there is a function $\limit^{\alpha,\varp} = \limit^{\Xi,\bar P,\alpha,\varp}\colon {}^{W_{\bar P}} Q_{\alpha,\varp}\to \Por$ such that, for any
    \begin{itemize}
        \item $(\Xi,\bar P)\in\calY$,
        
        \item $i^*<\omega$,
        
        \item $(\alpha_i,\varp_i)\in \theta\times (0,1)_\Q$,
        
        \item $\bar r^i = \{r^i_\ell \colon \ell\in W_{\bar P}\} \in {}^{W_{\bar P}}Q_{\alpha_i,\varp_i}$ for $i<i^*$,
        
        \item a finite partition $\la B_m \colon m<m^*\ra$ of $K_{\bar P}$,
        
        \item $\varp'>0$, and
        
        \item $q\in \Por$ stronger than $\limit^{\alpha_i,\varp_i} \bar r^i$ for all $i<i^*$,
    \end{itemize}
    there are some $q'\leq q$ in $\Por$, some non-empty finite $u\subseteq K_{\bar P}$ and some probability measure $\Xi^-$ on $\pts(u)$ such that
    \begin{enumerate}[label = \normalfont (\arabic*)]
        \item\label{i4II1} $|\Xi^-(u\cap B_m) - \Xi(B_m)| < \varp'$ for all $m<m^*$, and
        \item\label{i4II2} $\displaystyle \sum_{k\in u}\frac{|\{\ell\in I_k \colon q' \leq r^i_\ell\}|}{|I_k|} \Xi^-(\{k\}) > 1-\varp_i-\varp'$ for all $i<i^*$.
    \end{enumerate}
    \end{enumerate}
    Moreover, in~\ref{i4II}, whenever $\Xi$ is free, $\Xi^-$ can be the uniform measure on $u$.
\end{theorem}

\section{The proof of the main result}

In this section, we prove the main result \autoref{mt}.

For the rest of the section assume that: $K$ and $W$ are non-empty sets, $\bar{P} = \langle P_{k} \colon k \in K \rangle$ is a partition of $W$ into finite non-empty sets, and $\Xi$ is a finitely additive probability measure on $\calP(K).$  

\begin{lemma}\label{i18}   
    Let $\bbP$ be a forcing notion, $\langle B_{m} \colon m < m^{\ast} \rangle$  a finite partition of $K$, and $a_{m} \coloneqq \Xi(B_{m})$ for any $m < m^{\ast}.$ Let $r \in \bbP,$ $i^{\ast} < \omega,$ and $\langle \delta_{i} \colon i < i^{\ast} \rangle$ a sequence of real numbers. For any $i < i^{\ast}$ and $r' \in \bbP$  such that $r' \leq r,$ let $f_{r'}^{i} \colon K \to [0, 1]$ be a function. For $ m \in M \coloneqq \{ m < m^{\ast} \colon a_{m} > 0 \}$ and $ i< i^{\ast},$ define $$ c_{i, m}(r') \coloneqq \int_{K} f_{r'}^{i} d \Xi_{B_m} = \frac{1}{a_{m}} \int_{B_{m}} f_{r'}^{i} d \Xi.$$ If for any $r' \leq r$ and $i < i^{\ast} $ we have that $\int_{K} f_{r'}^{i} d \Xi \geq \delta_{i},$ then, for all $\varp > 0,$ there are some $r^{\ast} \leq r$ and a sequence of real numbers $\bar{c} = \langle c_{i, m} \colon i < i^{\ast} \conj  m \in M \rangle$ such that, for all $i < i^{\ast}$ and $m \in M$,
    
        \begin{enumerate}[label=\rm{(\arabic*)}]
            \item $0 \leq c_{i, m} \leq 1,$
    
            \item $\sum_{m \in M} c_{i, m} \cdot a_{m} \geq \delta_{i}$,
    
            \item $D^{\ast} \coloneqq  \left \{ r' \in \bbP \colon \forall i < i^{\ast} \, \forall m \in M \, ( \vert c_{i, m}(r') - c_{i, m} \vert < \varp) \right \}$ is dense below $r^{\ast}$ and $r^{\ast} \in D^{\ast}.$ 
        \end{enumerate}
\end{lemma}

\begin{PROOF}[\textbf{Proof}]{\autoref{i18}}
    Let $\varp > 0$ and $N < \omega$ large enough such that $\frac{1}{N} < \varp$ and $\calC \neq \emptyset,$ where $\calC$ is the set of finite sequences $\bar{c} = \langle c_{i, m} \colon i < i^{\ast} \conj  m \in M \rangle$ such that:

        \begin{itemize}
            \item $c_{i, m} \in [0, 1]_{\bbQ},$
    
            \item $N \cdot c_{i, m} \in \bbZ,$
    
            \item $\sum_{m \in M} c_{i, m} \cdot a_{m} \geq \delta_{i}.$ 
        \end{itemize}

    It is clear that $\calC$ is a finite set, hence there exists some $s^{\ast} < \omega$ such that $\calC = \{ \bar{c}^{s} \colon s < s^{\ast} \}.$ 

    Now, \textbf{suppose} that recursively on $s \leq s^{\ast}$, we can build a sequence $\langle r_{s} \colon s \leq s^{\ast} \rangle$ of conditions in $\bbP$  such that:
    \begin{itemize}
        \item $r_{0} \coloneqq r$ and for any $s < s^{\ast}, \, r_{s+1} \leq r_{s},$ 

        \item for any $r' \leq r_{s+1}, $ there are $i < i^{\ast}$ and $ m \in M $ such that $\vert c_{i, m}(r') - c_{i, m}^{s} \vert \geq \varp.$
    \end{itemize}
    For any $i<i^\ast$ and $m\in M$, since $c_{i,m}(r_{s^*})\leq 1$ and $\varp > \frac{1}{N}$, we can find some $c^\ast_{i,m}\in \left \{ \frac{\ell}{N} \colon  \ell \leq N \right\} \cap [ c_{i, m}(r_{s^{\ast}}),$ $c_{i, m}(r_{s^{\ast}}) + \varp)$. In particular, $\vert c_{i, m}^{\ast} - c_{i, m}(r_{s^{\ast}}) \vert < \varp$.

    
    Since $r_{s^{\ast}} \leq  r,$ by hypothesis and \autoref{t73}~\ref{t63.0}, we have $$\sum_{m \in M} c_{i, m}(r_{s^{\ast}}) \! \cdot \! a_{m} = \sum_{m \in M} \left( \frac{1}{a_{m}} \int_{B_{m}} f_{r_{s^{\ast}} }^{i} d \Xi \right) \! \cdot \! a_{m} \! = \! \sum_{m \in M} \left( \int_{B_{m}} f_{r_{s^{\ast}} }^{i} d \Xi \right) \! = \! \int_K f_{r_{s^{\ast}} }^{i} d \Xi \geq \delta_{i},$$
    
    hence, $$\delta_{i} \leq \sum_{m \in M} c_{i, m} (r_{s^{\ast}}) \cdot a_{m} \leq  \sum_{m \in M} c_{i, m}^{\ast} \cdot a_{m},$$ that is, $\bar{c}^{\ast} \coloneqq \{ c_{i,m}^{\ast} \colon i < i^{\ast} \conj  m \in M \} \in \calC,$ therefore there is some $s < s^{\ast}$ such that $\bar{c}^{\ast} = \bar{c}^{s}.$  
    Now, since $r_{s^{\ast}} \leq r_{s+1},$ we have that $\vert c_{i, m}^{s} - c_{i, m}(r_{s^{\ast}}) \vert \geq \varp$ for some $i < i^{\ast}$ and $m \in M.$ However, by the construction of $\bar{c}^{\ast},$ we have that  $ \vert c_{i, m}^{s} - c_{i, m}(r_{s^{\ast}}) \vert = \vert c_{i, m}^{\ast} - c_{i, m}(r_{s^{\ast}}) \vert < \varp,$ which is a contraction. Therefore, $r_{s^{\ast}}$ contradicts the choice of $r_{s+1}$ and, as a consequence, we cannot reach step $s^{\ast}$ in the induction, so we are ``stuck'' at some step $s < s^{\ast}.$
    
    Set $\bar c\coloneqq \bar c^s$. Let us show that $D^{\ast}$ is dense below $r_{s}$: let $q \leq r_{s}.$ Since  $r_{s+1}$ cannot be defined, there exists some $r' \leq q$ such that, for any $i < i^{\ast}$ and $m \in M,$ $\vert c_{i, m}^{s} - c_{i, m}(r') \vert < \varp,$ that is, $r' \in D^{\ast}.$ 
    
    Finally,  let $r^{\ast} \in D^{\ast}$ such that $r^{\ast} \leq r_{s}.$ It is clear that $r^{\ast},$ $\bar{c}$ are as required. 
\end{PROOF}

\autoref{i19} and~\ref{i23} below are the cornerstone to prove the main result (\autoref{mt}) of this paper. Although it is enough to restrict ourselves to Boolean algebras (i.e.\ $\Por=\cB$ and $\iota$ as the identity function), we develop our results for forcing notions densely embedding in a Boolean algebra with some properties, expecting that they can be useful to find new examples of $\theta$-$\FAM$-linked forcing notions (see more details in \autoref{i33}). However, for this paper, it does not hurt to restrict to $\Por=\cB$ and $\iota$ being the identity function.

The following result is a generalization of \cite[Lem.~2.17]{Sh00}. 

\begin{lemma}\label{i19}
     Assume that $\bbP$ is a forcing notion, $\cB$ is a Boolean algebra with a strictly positive finitely additive probability measure $\mu$, and $\iota \colon \bbP \to \cB^{+}$ is a dense embedding,  
     and $\Xi$ is a free finitely additive probability measure on $\calP(K).$ 
     Let $\langle B_{m} \colon m < m^{\ast} \rangle$ be a finite partition of $K,$ and for any $m < m^{\ast},$ $a_{m} \coloneqq \Xi(B_{m}).$ Let $r \in \bbP,$ $i^{\ast} < \omega,$ $\langle \delta_{i} \colon i < i^{\ast} \rangle$ a sequence of real numbers and for each $i < i^{\ast},$ $\bar{r}^{i} = \langle r_{\ell}^{i} \colon \ell \in W \rangle \in {}^{W} \bbP.$ If for any $ r' \leq r$ and $i < i^{\ast},$ $ \int_{K} f_{r'}^{i} d \Xi \geq \delta_{i},$ where $f_{r'}^{i} \colon K \to [0, 1]$ is defined by $f_{r'}^{i}(k) \coloneqq \frac{1}{\vert P_{k} \vert} \sum_{\ell \in P_{k}} \mu_{\iota(r')}(\iota(r_{\ell}^{i})),$ then, for all $\varp > 0$ and any finite set $F \subseteq K,$ there are a non-empty finite set $u \subseteq K \setminus F$ and some $r^{\oplus} \leq r$ such that: 

    \begin{enumerate}[label=\rm{(\arabic*)}]
        \item\label{i19.1}  $\left| \frac{\vert u \cap B_{m} \vert}{ \vert u \vert}  - \Xi(B_{m}) \right| < \varp,$ for all $m < m^{\ast}.$

        \item\label{i19.2} $\frac{1}{\vert u \vert} \sum_{k \in u} \frac{\vert \{ \ell \in P_{k} \colon  r^{\oplus} \leq r_{\ell}^{i} \} \vert}{\vert P_{k} \vert} \geq \delta_{i} - \varp,$  for all $i < i^{\ast}.$   
    \end{enumerate}
\end{lemma}

\begin{PROOF}[\textbf{Proof}]{\autoref{i19}}
    As in \autoref{i18}, we define $M \coloneqq \{ m < m^{\ast} \colon a_{m} > 0 \},$ hence $\sum_{m \in M} a_{m} = 1,$ because $\langle{B_{m} \colon m < m^{\ast}} \rangle$ is a partition of $W.$ Now, define for $m \in M,$ $i < i^{\ast} $ and $ r' \in \bbP,$ $ c_{i, m}(r') \coloneqq \frac{1}{a_{m}} \int_{B_{m}} f_{r'}^{i} d \Xi.$ Fix $\varp > 0$ and a finite set $F \subseteq K.$ Also consider, for any $r' \in \bbP$ and $i < i^{\ast},$ the map $ \varrho_{r'}^{i} \colon K \to \bbR$ such that $\varrho_{r'}^{i}(k) =  \frac{\vert \{ \ell \in P_{k} \colon r'  \leq r_{\ell}^{i} \} \vert}{\vert P_{k} \vert}.$ 

    Since for all $m \in M, \,  \frac{2 a_{m}(1 - a_{m})}{\varp^{2}} \geq 0,$ and $\frac{(\frac{\varp}{2})^{2}}{2(m^{\ast}+ i^{\ast})} > 0,$ there exists some $h^{\ast} < \omega$ such that $h^{\ast}$ is even and
    \begin{equation}\label{e149}
        \text{$\frac{2 \, a_{m}(1-a_{m})}{h^{\ast} \varp^{2}} < \frac{1}{m^{\ast} + i^{\ast}} \ \text{ and }  \ \frac{1}{h^{\ast}} < \frac{(\frac{\varp}{2})^{2}}{2(m^{\ast} + i^{\ast})}$}   
    \end{equation}

    On the other hand, since we choose $h^{\ast}$ such that $\frac{1}{h^{\ast}}  < \frac{(\frac{\varp}{2})^{2}}{2(m^{\ast}+ i^{\ast})},$ there exists some  $\varp^{\ast} > 0$ such that $\varp^{\ast} < \varp$ and:
    \begin{equation}\label{e150}
        \text{$ \frac{\frac{2}{h^{\ast}} + \varp^{\ast}}{(\frac{\varp}{2})^{2}} < \frac{1}{m^{\ast} + i^{\ast}} $}
    \end{equation}
    By \autoref{i18} applied to $\cB^{+},$   $\varp^{\ast}$ and the functions $f_{r'}^{i},$ there are $r^{\ast} \leq r$ and a sequence $\bar{c} = \langle c_{i, m} \colon i < i^{\ast} \conj  m \in M \rangle$ such that  $0 \leq c_{i, m} \leq 1,$  $\sum_{m \in M} c_{i, m} a_{m} \geq \delta_{i}$, the set $D^{\ast} \coloneqq  \left \{ r' \in \bbP \colon \forall i < i^{\ast} \, \forall m \in M \, ( \vert c_{i, m}(r') - c_{i, m} \vert <  \frac{\varp^{\ast}}{4}) \right \}$ is dense below $r^{\ast}$ and $r^{\ast} \in D^{\ast}.$ 
    
    The rest of the proof is too technical and quite long, so we will split it into several parts to make it easier to understand. First, we define a probability tree, and then instead of trying to find the desired $u$ and $r^{\oplus}$ directly, in the spirit of the probabilistic method, we measure the probability of their existence using the tree. Finally, we conclude that this probability is positive, so we can find such objects satisfying the needed requirements. We start with the construction of the tree. 

    \textbf{Part 1}: The tree construction.

    We set $\calF \coloneqq  \bigcup_{u \in [K]^{< \omega}} {}^{ \left( i^{\ast} \times \bigcup_{k \in u} P_{k} \right)} 2$ and $\pi_{2} \colon \calF \times K \to K$ such that $\pi_{2}(\sigma, k) \coloneqq k,$ that is, $\pi_{2}$ is the second component projection.

    Now, let us build, by recursion on the height $h$, a subtree $\calT$ of ${}^{< \omega}(M \cup (\calF \times K)),$ a function  $\bfr \colon \calT \to \bbP$  such that, for any $\rho \in \calT, \, \bfr(\rho) \coloneqq r_{\rho}$ and a probability space on $\suc(\rho)$ for any $\rho \in \calT.$  To illustrate the construction, see \autoref{f50}. Notice that, there,  $\eta = \langle m, (\sigma, k), m' \rangle.$
    
    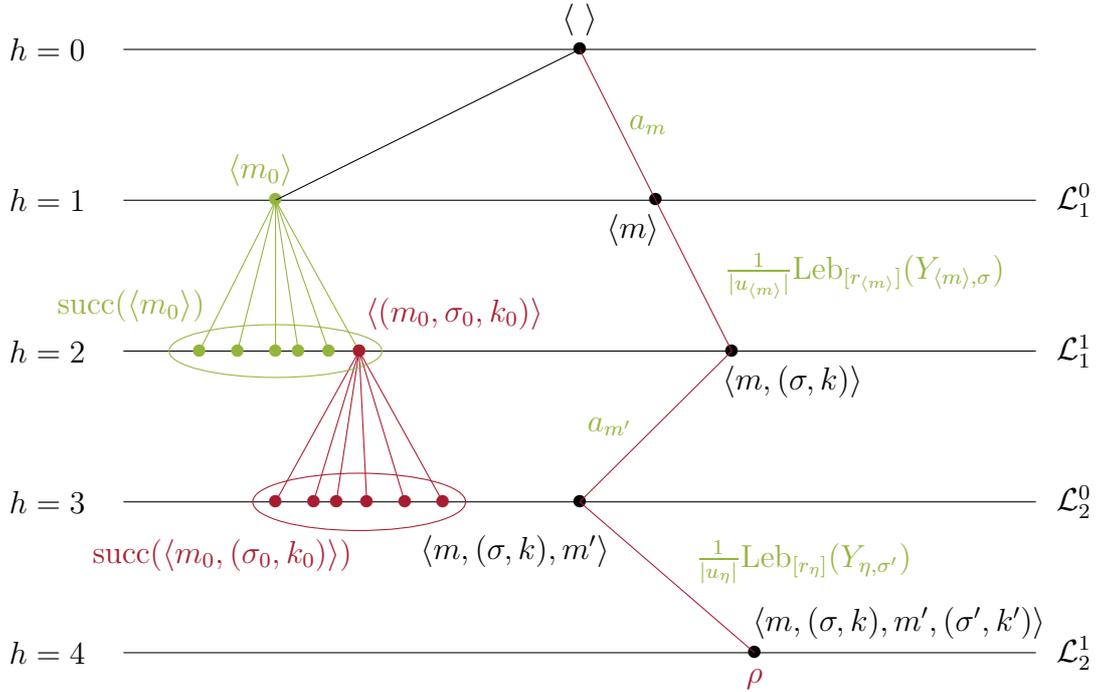
\begin{figure}[ht]
        \centering
        \begin{tikzpicture}
            \draw (0, 0) -- (12,0);
            \draw (0, -2) -- (12,-2);
            \draw (0, -4) -- (12,-4);
            \draw (0, -6) -- (12,-6);
            \draw (0, -8) -- (12,-8);
            \node at (6, 0.4) {$\langle \ \rangle$}; 
            \node at (-1,0) {$h = 0$}; 
            \node at (-1,-2) {$h = 1$}; 
            \node at (-1,-4) {$h = 2$}; 
            \node at (-1,-6) {$h = 3$}; 
            \node at (-1,-8) {$h = 4$}; 
            \node at (12.5,-2) {$\calL_{1}^{0}$}; 
            \node at (12.5,-4) {$\calL_{1}^{1}$}; 
            \node at (12.5,-6) {$\calL_{2}^{0}$}; 
            \node at (12.5,-8) {$\calL_{2}^{1}$}; 
            \node at (6,0) {$\bullet$}; 
            \node at (7,-2) {$\bullet$}; 
            \node at (8,-4) {$\bullet$}; 
            \node at (6,-6) {$\bullet$}; 
            \node at (8.3,-8) {$\bullet$};
        
            \draw[redun] (6, 0) -- (7,-2);
            \draw[redun] (7, -2) -- (8, -4);
            \draw[redun] (8, -4) -- (6,-6);
            \draw[redun] (6, -6) -- (8.3,-8);
        
            \node[greenun] at (6.9, -1) {$a_{m}$};
            \node at (6.7, -2.4) {$ \langle m \rangle$}; 
            \node[greenun] at (9.75, -3) {$\frac{1}{\vert u_{\langle m \rangle} \vert} \Leb_{[r_{\langle m \rangle}]}(Y_{\langle m \rangle, \sigma})$}; 
            \node at (8.8, -4.4) {$ \langle m, (\sigma, k) \rangle$}; 
        
            \node at (5.15, -6.65) {$\langle m, (\sigma, k), m' \rangle$};
            
            \node[greenun] at (6.4, -5) {$a_{m'}$};
            \node[redun] at (8.3, -8.35) {$\rho$};
            \node at (10.2, -7.6) {$\langle m, (\sigma, k), m', (\sigma', k') \rangle$};

            
            \node[greenun] at (2, -2) {$\bullet$};
            \draw (6, 0) -- (2, -2);

            \node[greenun] at (1, -4) {$\bullet$};
            \draw[greenun] (2, -2) -- (1, -4);
            \node[greenun] at (8.95, -6.8) {$\frac{1}{\vert u_{\eta} \vert} \Leb_{[r_{\eta}]}(Y_{\eta, \sigma'})$}; 
            \node[greenun] at (2, -4) {$\bullet$};
            \draw[greenun] (2, -2) -- (2, -4);

            \node[greenun] at (2.3, -4) {$\bullet$};
            \draw[greenun] (2, -2) -- (2.3, -4);

            \node[greenun] at (1.5, -4) {$\bullet$};
            \draw[greenun] (2, -2) -- (1.5, -4);

            \node[greenun] at (3.1, -4) {$\bullet$};
            \draw[greenun] (2, -2) -- (3.1, -4);

            \node[greenun] at (2.7, -4) {$\bullet$};
            \draw[greenun] (2, -2) -- (2.7, -4);
            \node[redun] at (3.1, -4) {$\bullet$};

            \node[redun] at (2, -6) {$\bullet$};
            \draw[redun] (3.1, -4) -- (2, -6);

            \node[redun] at (2.5, -6) {$\bullet$};
            \draw[redun] (3.1, -4) -- (2.5, -6);

            \node[redun] at (2.8, -6) {$\bullet$};
            \draw[redun] (3.1, -4) -- (2.8, -6);

            \node[redun] at (3.2, -6) {$\bullet$};
            \draw[redun] (3.1, -4) -- (3.2, -6);
            
            \node[redun] at (3.7, -6) {$\bullet$};
            \draw[redun] (3.1, -4) -- (3.7, -6);

            \node[redun] at (4.2, -6) {$\bullet$};
            \draw[redun] (3.1, -4) -- (4.2, -6);

            \draw[redun] (3.1,-6) ellipse (1.4cm and 0.38cm);

            \node[redun] at (1.3, -6.7) {$\suc(\langle m_{0}, (\sigma_{0}, k_{0}) \rangle)$};

            \node[redun] at (4.35, -3.5) {$\langle (m_{0}, \sigma_{0}, k_{0}) \rangle$};
            \draw[greenun] (2,-4) ellipse (1.4cm and 0.35cm);

            \node[greenun] at (0.1, -3.4) {$\suc(\langle m_{0} \rangle)$};

            \node[greenun] at (1.8, -1.6) {$\langle m_{0} \rangle$};
        \end{tikzpicture}
        \caption{A graphic example of the early levels of $\calT$.}
        \label{f50}
    \end{figure}    

    In the base step, we define $\Lev_{0}(\calT) \coloneqq \{ \langle \ \rangle \}$ and $r_{\langle \rangle} \coloneqq r^{\ast}.$ For the inductive step, suppose we have built the first $h$ levels and define the level $h + 1$ of $\calT$. For this, we will consider two possible cases:

    \textbf{Case 1}: When $h$ is even, $ \Lev_{h+1}(\calT) \coloneqq \calL_{\frac{h+2}{2}}^{0} \coloneqq \{ \rho^{\smallfrown} \langle m \rangle \colon m \in M \conj  \rho \in \Lev_{h}(\calT) \}.$ If $\eta \in \Lev_{h+1}(\calT),$ then there exists a unique $\rho \in \Lev_{h}(\calT)$ such that $\eta = \concat{\rho}{\langle m \rangle}.$ In this case, we define $ r_{\eta} \coloneqq r_{\rho}.$ In order to  define the probability space, let $\rho \in \Lev_{h}(\calT)$ and $\eta \in \suc(\rho),$ hence, there exists $m \in M$ such that $\eta = \concat{\rho}{\langle m \rangle}.$ In this case, we define $\p_{\rho}(\{ \eta \}) \coloneqq a_{m}.$ Thus, since $ \sum_{\eta \in \suc(\rho)} \p_{\rho}( \eta) = \sum_{m \in M} a_{m} = 1,$ it follows that $(\suc(\rho), \calP(\suc(\rho)), \p_{\rho})$ is a finite probability space.   

    \textbf{Case 2}: When $h$ is odd, we must work hard. 

    Let $\rho \in \Lev_{h}(\calT).$ Since $h -1$ is even, by the previous case, there exists $m \in M$ such that $\rho = \nu ^{\smallfrown} \langle m \rangle$ where $\nu \in \Lev_{h-1}(\calT).$ We set $$ \calK \coloneqq \{ \pi_{2}(\rho'(i)) \colon \rho' \in \Lev_{h}(\calT) \conj  0 \leq i < h \text{ is odd} \} \cup F.$$

    Since $\Xi$ is a free finitely additive measure and $\calK \subseteq K$ is a finite set, by \autoref{t90} there is a non-empty finite set $u_{\rho} \subseteq B_{m}$ such that $u_{\rho} \cap \calK = \emptyset,$ and
        \begin{equation}\label{152}
            \forall i < i^{\ast} \left( \left  \vert c_{i, m}(r_{\rho})- \frac{1}{ \vert u_{\rho} \vert} \sum_{k \in u_{\rho}} f_{r_{\rho}}^{i}(k)  \right \vert  < \frac{\varp^{\ast}}{4} \right)\! .
        \end{equation}
    Now, consider $\cB_{\rho}$ as the sub-Boolean algebra of $\cB_{\leq \iota(r_{\rho})}$ generated by the sets $ \{ \iota(r_{\rho}) \wedge \iota(r_{\ell}^{i}) \colon i < i^{\ast}, \ell \in P_{k}, k \in u_{\rho} \}.$ Defining $ y_{\rho, \sigma} \coloneqq \bigwedge_{(i, \ell) \in \dom(\sigma)} \iota(r_{\rho}) \wedge \iota(r_{\ell}^{i})^{\sigma(i, \ell)}$ for any $ \sigma \in {}^{i^{\ast} \times \bigcup_{k \in u_{\rho}} P_{k}} 2$ we have that, if  $y_{\rho, \sigma} \neq 0_\cB,$ then  
        \begin{multicols}{2}
            \begin{enumerate}
                \item[$\bullet$] $y_{\rho, \sigma} \leq_{\cB} \iota(r_{\ell}^{i}) \Leftrightarrow \sigma(i, \ell) = 0,$
    
                \item[$\bullet$] $y_{\rho, \sigma} \wedge_{\cB} \iota(r_{\ell}^{i}) = 0_\cB \Leftrightarrow \sigma(i, \ell) = 1.$
            \end{enumerate}
        \end{multicols}

        Also, since $\cB_{\rho}$ finitely generated and $\mu$ is strictly positive, by \autoref{b50}, we have that $\At_{\cB_{\rho}} = \{ y_{\rho, \sigma} \colon \sigma \in \Sigma_{\rho} \},$ where $\Sigma_{\rho} \coloneqq \{ \sigma \in 2^{i^{\ast} \times \bigcup_{k \in u_{\rho}} P_{k}} \colon \mu_{\iota(r_{\rho})}(y_{\rho, \sigma}) > 0  \}.$ As a consequence, 
        \begin{equation}\label{e154}
            \text{$\sum_{\sigma \in \Sigma_{\rho}} \mu_{\iota(r_{\rho})}(y_{\rho, \sigma}) = 1. $}
        \end{equation}
        
        So we get $u_{\rho}$ and $\Sigma_{\rho}$ for any $\rho \in \Lev_{h}(\calT),$ from which we can define: $$ \Lev_{h+1}(\calT) \coloneqq \calL_{\frac{h+1}{2}}^{1} \coloneqq \{ \rho^{\smallfrown} \langle (\sigma, k) \rangle \colon \rho \in \Lev_{h}(\calT), \, \sigma \in \Sigma_{\rho} \conj  k \in u_{\rho} \}. $$ 
        
        To define $\bfr$ at this level, let $\eta \in \Lev_{h+1}(\calT),$ hence there is some $\rho \in \Lev_{h}(T)$ such that $\eta = \rho{}^{\smallfrown} \langle (\sigma, k) \rangle$ for some $(\sigma, k) \in \Sigma_{\rho} \times u_{\rho}.$ By density, we can choose $r_{\eta} \in \bbP$ such that $\iota(r_{\eta}) \leq_{\cB} y_{\rho, \sigma}, \, r_{\eta} \leq r_{\rho}$ and $r_{\eta} \in D^{\ast}.$ Therefore, $\sigma(i, \ell) = 0$ implies $r_{\eta} \leq r_{\ell}^{i}$ and  $\sigma(i, \ell) = 1$ implies $ r_{\eta} \, \bot_{\bbP} \,  r_{\ell}^{i}.$

        In order to define the probability space, fix $\rho \in \Lev_{h}(\calT).$ For any $\concat{\rho}{\langle (\sigma, k) \rangle} \in \suc(\rho),$ define $\p_{\rho}( \concat{\rho}{\langle (\sigma, k) \rangle}) \coloneqq \frac{1}{\vert u_{\rho} \vert} \mu_{\iota(r_{\rho})}(y_{\rho, \sigma}).$ To prove that indeed $(\suc(\rho), \calP(\suc(\rho)), \p_{\rho})$ is a probability space, it is enough to show that $\sum_{\eta  \in  \suc(\rho)} \p_{\rho} (\eta) = 1.$ By \eqref{e154},  $$ \sum_{\eta \in \suc(\rho)}  \p_{\rho}( \eta)  = \sum_{\sigma \in \Sigma_{\rho}} \left( \mu_{\iota(r_{\rho})}(y_{\rho, \sigma}) \sum_{k \in u_{\rho}} \frac{1}{\vert u_{\rho} \vert} \right)\  = \sum_{\sigma \in \Sigma_{\rho}} \mu_{\iota(r_{\rho})}(y_{\rho, \sigma}) = 1.$$ Thus, $(\suc(\rho), \calP(\suc(\rho)), \p_{\rho})$ is a finite probability space. 

        Finally, we define $\calT \coloneqq \bigcup_{h < \omega} \Lev_{h}(\calT)$ and $\bfr \colon \calT \to \bbP$ such that $\bfr(\rho) \coloneqq r_{\rho}.$  It is clear by construction that $\calT$ is a probability tree.   

    For simplicity, we will use the following notation: for any $h \leq h^{\ast}$ consider  the sets $E_{h} \coloneqq \{ n < \omega \colon 1 \leq n \leq h \text{ even} \}  $ and $   O_{h} \coloneqq  \{ n < \omega \colon 1 \leq n \leq h \text{ is odd} \}.$ Also,  if $\rho \in \Lev_{h}(\calT)$ for $h \geq 2,$ then by construction either  $\rho = \eta ^{\smallfrown} \langle m, (\sigma, k) \rangle$ or  $\rho = \eta ^{\smallfrown} \langle (\sigma, k), m \rangle.$ In any case, we denote $ m_{\rho} \coloneqq m,  \sigma_{\rho} $ $\coloneqq \sigma$ and $ k_{\rho} \coloneqq k.$ That is, $m_{\rho}, \, k_{\rho},  \, \sigma_{\rho}$ are the last $m, k, \sigma$ that appear in $\rho.$ \index{$m_{\rho}$} For $\rho$ of even length $\geq 2,$ set $u_{\rho} \coloneqq u_{\rho \, {\rest} \, h-1},$ and allow $m_{\rho} = m$ when $\rho = \langle m \rangle.$  \index{$k_{\rho}$} \index{$\sigma_{\rho}$} As a consequence, for $k \in u_{\rho}$ where $\rho$ of even length and $i < i^{\ast},$ $$  \varrho_{\rho}^{i}(k) \coloneqq \varrho_{r_{\rho}}^{i}(k) =  \frac{\vert \{ \ell \in P_{k} \colon r_{\rho}  \leq r_{\ell}^{i} \} \vert}{\vert P_{k} \vert} =   \frac{\vert \{ \ell \in P_{k} \colon \sigma_{\rho}(i, \ell) = 0 \} \vert}{\vert P_{k} \vert}.$$

    Since $\calT$ is a probability tree, by  \autoref{p49}, it induces a probability space in each of its levels.

    \textbf{Part 2}: Finding a suitable $\rho \in \Lev_{h^{\ast}}(\calT)$ with high probability.

    In general, as the length of the elements of $\calT$ increases probabilities are getting smaller and smaller, what we do is try to define random variables that allow us to model what happens in both condition~\ref{i19.1} and condition~\ref{i19.2} of the conclusion of the lemma. For this, we divide this part into two sub-parts: in the first, we deal with condition~\ref{i19.1}, and in the second, with condition~\ref{i19.2}.
    
    For $\rho \in T,$ define $u_{\rho}^{\ast} \coloneqq \{ k_{\rho \rest h} \colon h \in E_{h^{\ast}} \}$.\footnote{We are counting only the even numbers because the information is repeated in the construction, for example, if $h$ is even, at the levels $h$ and $h+1$ we have the same $k_{\rho}$ for any $\rho \in \Lev_{h}(\calT) \cup \Lev_{h+1}(\calT).$ So adding over all numbers below $h^{\ast}$ would affect the value of the probabilities.} It is clear by construction that $\vert u_{\rho}^{\ast} \vert = \frac{h^{\ast}}{2}, \, u_{\rho}^{\ast} \subseteq \bigcup_{m \in M} B_{m}$ and $u_{\rho}^{\ast} \subseteq K \setminus F.$  

    \textbf{Part 2.1}: Random variables to model the first condition.
    
    For $m \in m^{\ast} \setminus M, \, u_{\rho}^{\ast} \cap B_{m} = \emptyset,$  therefore $ \left \vert \frac{ \vert u_{\rho}^{\ast} \cap B_{m} \vert}{ \vert u_{\rho}^{\ast} \vert} - a_{m} \right \vert = 0 < \varp.$  

    On the other hand, for $m \in M$  consider $U_{m}, V_{m} \colon \Lev_{h^{\ast}}(\calT) \to \bbR$ such that, if $\rho \in  \Lev_{h^{\ast}}(\calT),$ then $V_{m}(\rho) \coloneqq \frac{\vert u^{\ast}_{\rho} \cap B_{m} \vert}{\vert u_{\rho}^{\ast} \vert} \ \text{ and } \ U_{m}(\rho) \coloneqq \vert \{ j \in E_{h^{\ast}} \colon m_{\rho \rest j} = m \}  \vert.$ Since the $\sigma$-algebra in $\Lev_{h^{\ast}}(\calT)$ is $\calP(\Lev_{h^{\ast}}(\calT)),$ $U_{m}$ and $V_{m}$ are, trivially, random variables. Also, $U_{m} \sim \Bin(\frac{h^{\ast}}{2}, a_{m}),$ since it counts the success after $\frac{h^{\ast}}{2}$ tries, each with probability $a_{m},$ and we can express $V_{m}$ in terms of $U_{m}$:  $V_{m}(\rho) = \frac{2 U_{m}(\rho)}{h^{\ast}}.$ Hence, $$ \E[V_{m}] = \frac{2 \, \E[U_{m}]}{h^{\ast}} = \frac{ \frac{2 \, h^{\ast}}{2} \, a_{m}}{h^{\ast}} = a_{m} \text{ and } \Var[V_{m}] = \frac{\frac{4 \, h^{\ast} }{2}\, a_{m} \, (1 - a_{m})}{(h^{\ast})^{2}} = \frac{2 \, a_{m} (1 - a_{m})}{h^{\ast}}.$$ 
    
    Thereby, by \emph{Chebyshev inequality}, $\p_{h^{\ast}}[\vert V_{m} - a_{m} \vert \geq \varp] \leq \frac{\Var[V_{m}]}{\varp^{2}} = \frac{ 2 \, a_{m} \, (1 - a_{m})}{h^{\ast} \, \varp^{2}}.$ Thus, by the choice of $h^{\ast},$\footnote{In fact, \eqref{e156} was what motivated the choice of $h^{\ast}$ in \eqref{e149}.} we get:  
    \begin{equation}\label{e156}
        \text{$\forall m \in M \left( \p_{h^{\ast}}[\vert V_{m} - a_{m} \vert \geq \varp]  \leq \frac{1}{m^{\ast} + i^{\ast}} \right)\! .$}
    \end{equation}

    In the next part, we will do the same for the second condition.

    \textbf{Part 2.2}: Random variables to model the second condition.

    Now, for $h \in E_{h}$ and $i < i^{\ast},$ consider $Z_{h}^{i} \colon \Lev_{h}(\calT) \to \bbR$ such that $Z_{h}^{i}(\eta) \coloneqq \varrho_{\eta}^{i}(k_{\rho})$ for any $\eta \in \Lev_{h}(\calT).$ For $h \leq  h^{\ast}$ even and $\rho \in \Lev_{h-1}(\calT),$ we calculate the relative expected value of $Z_{h}^{i}$. For this,  since $h-1$ is odd, the successors of $\rho$ are as in \textbf{Case 2} of Part 1 and we use  \autoref{p60}: 
    \vspace{-0.4cm}
    \begin{equation*}
        \begin{split}
            \E[Z_{h}^{i} \colon  \eta \, {\rest}_{h-1} = \rho] & = \E_{\suc(\rho)}[Z_{h}^{i} \rest \suc(\rho)] = \sum_{k \in u_{\rho}} \left( \sum_{\sigma \in \Sigma_{\rho}} \frac{1}{\vert u_{\rho} \vert} \mu_{\iota(r_{\rho})} (y_{\rho, \sigma}) \,  Z_{h}^{i}(\concat{\rho}{\langle (\sigma, k) \rangle}) \right)\\
            & = \sum_{k \in u_{\rho}} \left( \sum_{\sigma \in \Sigma_{\rho}} \frac{1}{\vert u_{\rho} \vert}{\mu_{[r_{\rho}]}(Y_{\rho, \sigma})} \, \varrho_{\rho^{\smallfrown} \langle  (\sigma, k) \rangle}^{i} (k) \right)\\
            & = \frac{1}{\vert u_{\rho} \vert} \sum_{k \in u_{\rho}} \left( \sum_{\sigma \in \Sigma_{\rho}} \mu_{\iota(r_{\rho})}(y_{\rho, \sigma}) \,   \frac{ \vert \{ \ell \in P_{k} \colon \sigma(i, \ell) = 0 \} \vert}{\vert P_{k} \vert}  \right)  \\
         & = \frac{1}{\vert u_{\rho} \vert} \sum_{k \in u_{\rho}} \left[ \frac{1}{\vert P_{k} \vert} \sum_{\sigma \in \Sigma_{\rho}} \mu_{\iota(r_{\rho})}(y_{\rho, \sigma}) \left( \sum_{\ell \in P_{k}, \ \sigma(i, \ell) = 0} 1 \right)  \right]\\
            & = \frac{1}{\vert u_{\rho} \vert} \sum_{k \in u_{\rho}} \left[ \frac{1}{\vert P_{k} \vert} \sum_{\ell \in P_{k}} \left(  \sum_{\sigma \in \Sigma_{\rho}, \ \sigma(i, \ell) = 0} \mu_{\iota(r_{\rho})} (y_{\rho, \sigma})  \right)  \right]\\
            & = \frac{1}{\vert u_{\rho} \vert} \sum_{k \in u_{\rho}} \left[ \frac{1}{\vert P_{k} \vert} \sum_{\ell \in P_{k}}  \mu_{\iota(r_{\rho})}(\iota(r_{\ell}^{i}))  \right] = \frac{1}{\vert u_{\rho} \vert} \sum_{k \in u_{\rho}} f_{r_{\rho}}^{i}(k).
        \end{split}
    \end{equation*}

    Therefore, by the choice of $u_{\rho}$ (see \eqref{152}),  $\left \vert \E  [ Z_{h}^{i} \colon  \eta \, {\rest}_{h-1} = \concat{\rho}{\langle m \rangle} ] - c_{i, m}(r_{\rho}) \right \vert <\frac{\varp^{\ast}}{4},$ for any $i < i^{\ast}$. As a consequence, since $r_{\rho} \in D^{\ast},$ we can conclude:
    \begin{equation}\label{e158}
        \text{$\forall i < i^{\ast} \left(  \left \vert \E  [ Z_{h}^{i} \colon  \eta \, {\rest}_{h-1} = \concat{\rho}{\langle m \rangle} ] - c_{i, m}^{\ast} \right \vert <\frac{\varp^{\ast}}{2} \right).$}
    \end{equation}

    Now, we are going to calculate the relative expected value restricted to one more level. Since $h - 2$ is even, the successors of $\nu \in \Lev_{h-2}(\calT)$ are as in the first case of Part 1. By \autoref{p64}, we have that: 
    \begin{equation*}
        \begin{split}
            \E[Z_{h}^{i} \colon  \eta \, {\rest}_{h-2} = \nu \, ] & = \E[  \E [Z_{h}^{i} \colon  \eta \, {\rest}_{h-1} = \rho] \colon  \rho \, {\rest}_{h-2} = \nu] = \E \left[ \frac{1}{\vert u_{\rho} \vert} \sum_{k \in u_{\rho}} f_{r_{\nu}}^{i}(k) \colon  \rho \, {\rest}_{h-2} = \nu \right]\\
            & = \E_{\suc(\nu)} \left[ \frac{1}{\vert u_{\rho} \vert} \sum_{k \in u_{\rho}} f_{r_{\nu}}^{i}(k)   \right] = \sum_{m \in M} \left( a_{m} \cdot \frac{1}{\vert u_{ \concat{\nu}{\langle m \rangle}} \vert} \sum_{k \in u_{\concat{\nu}{\langle m \rangle}}} f_{r_{\nu}}^{i}(k) \right)\!,
        \end{split}
    \end{equation*}
    and therefore,
    \begin{equation*}
        \begin{split}
            \left \vert \E[Z_{h}^{i} \colon  \eta \, {\rest}_{h-2} = \nu \, ] - \! \!   \int_{K} \! \!  f_{r_{\nu}}^{i} d  \Xi \right \vert & =    \left \vert  \sum_{m \in M}  \left( \frac{a_{m}}{\vert u_{ \concat{\nu}{\langle m \rangle}} \vert} \sum_{k \in u_{\concat{\nu}{\langle m \rangle}}}  f_{r_{ \nu}}^{i}(k) \right)  - \int_{K}  f_{r_{ \nu}}^{i} d \Xi \right \vert\\
            & \leq   \sum_{m \in M} a_{m} \, \left \vert \left( \frac{1}{\vert u_{ \concat{\nu}{\langle m \rangle}} \vert} \sum_{k \in u_{ \concat{\nu}{\langle m \rangle}}}  f_{r_{ \nu}}^{i}(k) - c_{i, m}(r_{ \concat{\nu}{\langle m \rangle}}) \right)  \right \vert\\
            & < \frac{\varp^{\ast}}{4}.
        \end{split}
    \end{equation*}

    Since $r_{\nu} \leq r,$ we have that  $\int_{K} f_{r_{\nu}}^{i}d\Xi \geq \delta_{i}$, hence, by \eqref{e158},
    \begin{equation}\label{e160}
        \text{$\forall i < i^{\ast} \left( \E[Z_{h}^{i} \colon  \eta \, {\rest}_{h-2} = \nu \, ] > \delta_{i} - \frac{\varp^{\ast}}{4} 
         \right) \!. $ }
    \end{equation}
    Similarly, we can prove that 
    \begin{equation}\label{e162}
        \text{$ \forall i < i^{\ast} \left(  \left \vert \E[Z_{h}^{i} \colon  \eta \, {\rest}_{h-2} = \rho \, ] - \sum_{m \in M} a_{m} c_{i, m}  \right \vert < \frac{\varp^{\ast}}{2}  \right)\!  .$}
    \end{equation}

    The previous equations can be generalized in the following way. For $ i < i^{\ast}, $ $ h \in  E_{h^{\ast}}, $ $j \in E_{h -1},$ and $\rho \in \Lev_{j}(\calT),$ we have that $\E[Z_{h}^{i} \colon  \eta \, {\rest}_{j} = \rho \, ] = \E_{h-2}[\E[Z^i_h \colon \eta\, {\rest}_{h-2} = \rho'] \colon \rho'\, {\rest}_{j} = \rho]$ by \autoref{p66}. On the other hand, $\Lev_{h-2}(\calT)$ is finite, so we can find some $\varp'>0$ such that the absolute value in \eqref{e162} is even smaller than $\frac{\varp^*}{2}-\varp'$ for any $\rho'\in \Lev_{h-2}(\calT)$, that is, $\sum_{m\in M}a_m c_{i,m} -\frac{\varp^*}{2} + \varp'  <  \E[Z_{h}^{i} \colon  \eta \, {\rest}_{h-2} = \rho' \, ] < \sum_{m\in M}a_m c_{i,m} +\frac{\varp^*}{2} - \varp'$. Hence, after taking the expected value on $\rho'$, $\sum_{m\in M}a_m c_{i,m} -\frac{\varp^*}{2} + \varp'  \leq  \E[Z_{h}^{i} \colon  \eta \, {\rest}_{j} = \rho \, ] \leq \sum_{m\in M}a_m c_{i,m} +\frac{\varp^*}{2} - \varp'$, that is, 
    %
    \begin{equation}\label{e163}
        \text{$  \left \vert \E[Z_{h}^{i} \colon  \eta \, {\rest}_{j} = \rho \, ] - \sum_{m \in M} a_{m} c_{i, m}  \right \vert < \frac{\varp^{\ast}}{2} $.}
    \end{equation}
     Similarly, by \eqref{e160} we obtain that $E[Z^i_h\colon \eta\, {\rest}_j = \rho] > \delta_i - \frac{\varp^*}{4}$. Thus, for $j=0$, we get:
    \begin{equation}\label{e164}
        \text{$\forall i < i^{\ast} \left( \E[Z_{h}^{i}] > \delta_{i} - \frac{\varp^{\ast}}{4}  \text{ and } \left \vert \E[Z_{h}^{i}] - \sum_{m \in M} a_{m} \, c_{i, m}  \right \vert < \frac{\varp^{\ast}}{2} \right) \!.$} 
    \end{equation}

    Now, for any $i < i^{\ast}$ consider the random variable $Y_{i} \colon \Lev_{h^{\ast}}(\calT) \to \bbR$ such that, for every $\rho \in \Lev_{h^{\ast}}(\calT),$ $ Y_{i} (\rho) \coloneqq \frac{1}{\vert u_{\rho}^{\ast} \vert} \sum_{k \in u_{\rho}^{\ast}} \frac{\vert \{ \ell \in P_{k} \colon r_{\rho} \leq r_{\ell}^{i} \}}{\vert P_{k} \vert}.$ Notice that we can express $Y_{i}$ in terms of $Z_{h}^{i}$ for $h \in E_{h^{\ast}}$ as follows: $ Y_{i} = \frac{2}{h^{\ast}} \sum_{h \in E_{h^{\ast}}} Z_{h}^{i}$ hence, by basic properties of the expected value and \eqref{e164}, we get: 
    \begin{equation*}
        \text{$\E[Y_{i}] = \E \left[ \frac{2}{h^{\ast}} \sum_{h \in E_{h^{\ast}}} Z_{h}^{i} \right] = \frac{2}{h^{\ast}} \sum_{h \in E_{h^{\ast}}} \E[Z_{h}^{i}]  > \frac{2}{h^{\ast}} \sum_{h \in E_{h^{\ast}}} \left( \delta_{i} - \frac{\varp^{\ast}}{4} \right) =  \delta_{i} - \frac{\varp^{\ast}}{4}.$}
    \end{equation*}
    
    On the other hand, since $\varp^{\ast} < \varp,$ it is clear that $\E[Z_{h}^{i}] > \delta_{i} - \frac{\varp}{4}$ and $\delta_{i} - \frac{\varp}{4} < \E[Y_{i}],$ hence: $$ \p[Y_{i} \leq \delta_{i} - \varp] \leq \p \left [Y_{i} \leq \E [ Y_{i}] - \frac{\varp}{2} \right] \leq \p \left [\frac{\varp}{2} \leq \vert Y_{i} - \E[Y_{i}] \vert \right],$$ and applying \emph{Chebyshev's inequality}, we get: 
    
    \begin{equation}\label{e165}
        \text{$\forall i <i^{\ast} \left( \p[Y_{i} \leq \delta_{i} - \varp] \leq \frac{ \Var[Y_{i}]}{(\frac{\varp}{2})^{2}} \right) \! .$}
    \end{equation}

    \textbf{Part 2.2.1}: Properly bound the variance of $Y_{i}.$ 

    Since we want $\p[ \forall i < i^{\ast}(Y_{i} > \delta_{i} - \varp)] > 0,$ we must show that $\Var(Y_{i})$ is small enough. For this, let us start by noting that, by \autoref{p39}~\ref{p39.3}:
    $$ \Var[Y_{i}] = \frac{4}{(h^{\ast})^{2}} \left( \sum_{h \in E_{h^{\ast}}} \Var[Z_{h}^{i}] + \sum_{ h, j \in E_{h^{\ast}}, \, h \neq j} \Cov[Z_{j}^{i}, Z_{h}^{i}] \right) \! .$$

    It is clear that, for any  $i < i^{\ast}, \, 0 \leq \Var[Z_{h}^{i}] \leq 1,$ because  $0 \leq Z_{h}^{i} \leq 1;$ and $0 \leq \E[Z_{h}^{i}] \leq 1$ and $\vert Z_{h}^{i} - \E[Z_{h}^{i}] \vert \leq 1.$  Therefore, we must bound the covariance. For this, let $j, h \in E_{h^{\ast}}$ such that $j < h.$ Then, by \autoref{p66}, \eqref{e164} and \eqref{e162},  we get:\footnote{Recall that the subindex in $\E$ refers to the level where the expected value is calculated. Notice that $\E_h[Z^i_j] = \E_j[Z^i_j]$ since the random variable $Z^i_j$ only depends on level $j$.}
    \begin{equation*}
        \begin{split}
            \Cov[Z_{j}^{i}, Z_{h}^{i}]& = \E_h[ Z_{j}^{i} \cdot Z_{h}^{i}] - \E_j[Z_{j}^{i}] \cdot \E_h[Z_{h}^{i}] = \E_{j} [ \, \E_h \, [ Z_{j}^{i}  \cdot Z_{h}^{i} \colon   \nu \rest_{j} \, = \, \eta \, ]]  - \E_{j}[Z_{j}^{i} \cdot \E_h \, [Z_{h}^{i}]]\\
            & = \E_{j} \left[  Z_{j}^{i} \cdot \left( \,  \E \,[ Z_{h}^{i}  \colon   \nu \rest_{j} \, = \, \eta ] - \E_{h}[Z_{h}^{i}]  \right)  \right]\\
            & \leq \E_{j} \left[  Z_{j}^{i} \cdot \left( \left \vert \,  \E[ Z_{h}^{i}   \colon   \nu \rest_{j} \, = \, \eta ] - \E_{h}[Z_{h}^{i}] \, \right \vert \right)  \right]\\
            & \leq \E_{j} \left[  Z_{j}^{i} \cdot \left( \left \vert \E[ Z_{h}^{i} \colon   \nu \rest_{j} \, = \, \eta ] - \sum_{m \in M} c_{i, m} \, a_{m} \right \vert + \left \vert E[Z_{h}^{i}] - \sum_{m \in M} c_{i, m} \, a_{m}  \right \vert \right)  \right]\\
            & \leq \E_{j} \left[ Z_{j}^{i} \cdot \left(  \frac{\varp^{\ast}}{2} + \frac{\varp^{\ast}}{2} \right) \right]  = \varp^{\ast} \cdot \E_{j}
            [Z_{j}^{i}] \leq \varp^{\ast}. 
        \end{split}
    \end{equation*}

    As a consequence, $$ \Var[Y_{i}] \leq \frac{4}{(h^{\ast})^{2}} \left[ \frac{h^{\ast}}{2} + \frac{h^{\ast}}{2} \left( \frac{h^{\ast}}{2}-1 \right)  \varp^{\ast} \right] = \frac{2}{h^{\ast}} + \varp^{\ast} - \frac{2}{h^{\ast}} \varp^{\ast} =  \frac{2}{h^{\ast}} + \left( \frac{h^{\ast} -2}{h^{\ast}} \right) \varp^{\ast} < \frac{2}{h^{\ast}} + \varp^{\ast}.$$ 
    
    Thus, \eqref{e150}, \eqref{e165} and the choice of $\varp^{\ast}$ entail:  
    \begin{equation}\label{e167}
        \text{$\forall i < i^{\ast} \left( \p[Y_{h^{\ast}}^{i} < \delta_{i} - \varp] < \frac{1}{m^{\ast} + i^{\ast}} \right) \!.  $}
    \end{equation}
    
    \textbf{Part 2.3}: Some $\rho$ of high probability and the conclusion.

    Consider the following events: 
    \begin{itemize}
        \item $E \coloneqq \{ \rho \in \Lev_{h^{\ast}}(\calT) \colon \forall m \in M \, \forall i < i^{\ast}  \left( \vert V_{m, h^{\ast}}(\rho) - a_{m} \vert < \varp \conj  Y_{i}(\rho) > \delta_{i} - \varp \right) \},$ 
        
        \item $F \coloneqq \{ \rho \in \Lev_{h^{\ast}}(\calT) \colon \exists m \in M ( \vert V_{m, h^{\ast}}(\rho)) - a_{m} \vert \geq \varp \},$

        \item $G \coloneqq \{ \rho \in \Lev_{h^{\ast}}(\calT) \colon \exists i < i^{\ast} (Y_{i}(\rho) \leq \delta_{i} - \varp \}.$
    \end{itemize}

    It is clear that $E = F^{\rm{c}} \cap G^{\rm{c}} = (F \cup G)^{\rm{c}}$ hence, $E^{\rm{c}} = F \cup G.$ Also, by \eqref{e156} and \eqref{e167} we  have that 
    $$ \p(F) \leq \sum_{m \in M} \p[\vert V_{m, h^{\ast}} - a_{m} \vert \geq \varp] \leq \sum_{m \in M} \frac{1}{m^{\ast} + i^{\ast}} \leq \frac{m^{\ast}}{m^{\ast} + i^{\ast}}, \text{ and } $$ $$ \p(G) \leq \sum_{i < i^{\ast}} \p[Y_{h^{\ast}}^{i} \leq \delta_{i} - \varp] < \sum_{i < i^{\ast}} \frac{1}{m^{\ast}+ i^{\ast}} = \frac{i^{\ast}}{m^{\ast} +i^{\ast}}.$$ 
    Therefore,
    \begin{equation*}
        \begin{split}
            \p(E) & = 1 - \p(E^{\rm{c}}) = 1 - [\p(F) + \p(G)] - \p(F \cap G)]\\
            & > 1 - \left( \frac{m^{\ast}}{m^{\ast}+ i^{\ast}} + \frac{i^{\ast}}{m^{\ast}+ i^{\ast}} \right) + \p(F \cap G)\\
            & = 1 - \frac{m^{\ast} + i^{\ast}}{m^{\ast} + i^{\ast}} + \p(F \cap G)  = \p(F \cap G) \geq 0,
        \end{split}
    \end{equation*}
    
    hence $\p(E) > 0.$ As a consequence, $E \neq \emptyset.$ Let $\rho \in E, \, u \coloneqq u_{\rho}^{\ast}$ and $r^{\oplus} \coloneqq \bfr(\rho) = r_{\rho}.$ Then, by the construction of $\calT$, $u \subseteq K \setminus F$ and $r^{\oplus} = r_{\rho} \leq r^{\ast} \leq r.$ Also, since $\rho \in E,$ for $m \in M$ and $i < i^{\ast},$ we have that:  
    $$ \left \vert \frac{\vert u \cap B_{m} \vert}{ \vert u \vert} - a_{m} \right \vert = \left \vert \frac{ \vert u_{\rho}^{\ast} \cap B_{m} \vert}{\vert u_{\rho}^{\ast} \vert} - a_{m} \right \vert = \vert V_{m, h^{\ast}}(\rho) - a_{m} \vert < \varp, \text{ and }$$ 
    $$ \frac{1}{\vert u \vert} \sum_{k \in u} \frac{\vert \{ \ell \in P_{k} \colon r^{\oplus} \leq r_{\ell}^{i} \}  \vert}{\vert P_{k} \vert} = \frac{1}{\vert u_{\rho }^{\ast} \vert} \sum_{k \in u_{\rho}^{\ast}} \frac{\{ \ell \in P_{k} \colon r_{\rho} \leq r_{\ell}^{i} \}}{\vert P_{k} \vert} = Y_{i}(\rho) > \delta_{i} - \varp.$$
    
    Thus, finally:   
    \begin{enumerate}[label=\rm{(\arabic*)}]
        \item  $\left \vert \frac{\vert u \cap B_{m} \vert}{ \vert u \vert}  - a_{m}) \right \vert < \varp,$ for all $m < m^{\ast}.$

        \item $\frac{1}{\vert u \vert} \sum_{k \in u} \frac{\vert \{ \ell \in P_{k} \colon  r^{\oplus} \leq r_{\ell}^{i} \} \vert}{\vert P_{k} \vert} > \delta_{i} - \varp,$ for all $i < i^{\ast}.$  \qedhere  
    \end{enumerate}
\end{PROOF}

The following result is a generalization of \cite[Lem.~2.18]{Sh00} and will allow defining the limit as in \autoref{i2}~\ref{i2II} in the general case (see the proof of \autoref{i27}). 

\begin{lemma}\label{i23}
    Assume that $\bbP$ forcing notion, $\cB$ is a $\sigma$-complete Boolean algebra with a strictly positive measure $\mu$,\footnote{Notice that this implies that $\cB$ is complete.} and $\iota \colon \bbP \to \cB^{+}$ is a dense embedding. Let $ \delta \in [0, 1],$ $r^{\ast} \in \bbP$ and $\bar{r} = \langle r_{\ell} \colon \ell \in W \rangle \in {}^{W}\bbP$ be such that, for any $\ell \in W,$ $ \mu_{\iota(r^{\ast})}(\iota(r_{\ell})) \geq \delta.$ Then there exists $r^{\otimes} \in \bbP$ such that $r^{\otimes}  \leq r^{\ast}$ and, for any $r \leq r^{\otimes},$ $  \int_{K} f_{\iota(r)}^{\iota(\bar{r})}  d \Xi \geq \delta,$ where $\iota(\bar{r}) \coloneqq \langle \iota(r_{\ell}) \colon \ell \in W \rangle$, and for any $b \in \cB$ and $\bar{b} = \langle b_{\ell} \colon \ell  \in W \rangle \in {}^{W} \cB$, $f_{b}^{\bar{b}} \colon K \to \bbR$ is defined by $ f_{b}^{\bar{b}}(k) \coloneqq \frac{1}{\vert P_{k} \vert} \sum_{\ell \in P_{k}} \mu_{b}(b_{\ell}).$ 
\end{lemma}

\begin{PROOF}[\textbf{Proof}]{\autoref{i23}}
    We consider the set $I \coloneqq \left \{  r \in \bbP \colon r \leq r^{\ast}  \conj  \int_{K} f_{\iota(r)}^{\iota(\bar{r})} d \Xi  < \delta \right  \}.$ Suppose that $I$ is not dense below $r^{\ast}.$ So there exists $r^{\otimes} \leq r^{\ast}$ such that, for all $r \in I,$ $r \nleq r^{\otimes},$ which implies that, for any $r \leq r^{\otimes},$ $r \notin I,$ that is, $r^{\otimes}$ is as required. So it is enough to prove that $I$ is not dense below $r^{\ast}.$ Towards contradiction, assume that $I$ is dense below $r^{\ast},$ hence we can find a maximal anti-chain $A \coloneqq \{ s_{i} \colon i < i_{\ast} \} \subseteq I$ below $r^{\ast}.$ Notice that $\bbP$ is c.c.c., hence we have that $0 < i_{\ast}  \leq \omega.$ Also, we can think of $A$ almost like a partition in the sense that $i, j < i_{\ast}$ and $i \neq j$ imply $\mu(\iota(s_{i}) \wedge \iota(s_{j})) = 0.$ As a consequence, $\mu(\iota(r^{\ast})) = \sum_{i < i_{\ast}} \mu(\iota(s_{i})).$ 
    
    Now, for all $j < \min \{ i^{\ast}+ 1, \omega \}$ we define $s^{j} \coloneqq \bigvee_{i < j}\iota(s_{i}),$ hence we can express $f_{\iota(s^{j})}^{\iota(\bar{r})}(k)$ in terms of $f_{\iota(s_{i})}^{\iota(\bar{r})}(k)$ for $k \in K,$ $0 < j < i_{\ast}$ and $i < j,$ as follows: $$ f_{\iota(s^{j})}^{\iota(\bar{r})}(k) = \sum_{i < j} \mu_{\iota(s^{j})}(\iota(s_{i})) \cdot f_{\iota(s_{i})}^{\iota(\bar{r})}(k).$$
    
    Since $s_{0} \in I,$ $\int_{K} f_{\iota(s_{0})}^{\iota(\bar{r})} d \Xi < \delta$ and $\varp \coloneqq \delta - \int_{K} f_{\iota(s_{0})}^{\iota(\bar{r})}  d \Xi > 0.$ Then, 
    \begin{equation*}
        \begin{split}
            \int_{K} f_{\iota(s^{j})}^{\iota(\bar{r})} d \Xi &= \int_{K} \left(  \sum_{i < j} \mu_{\iota(s^{j})}(\iota(s_{i})) f_{\iota(s_{i})}^{\iota(\bar{r})} \right) d \Xi  = \sum_{i < j} \mu_{\iota(s^{j})}(\iota(s_{i})) \int_{K} f_{\iota(s_{i})}^{\iota(\bar{r})} d \Xi \\
            & =  \mu_{\iota(s^{j})}(\iota(s_{0})) \int_{K} f_{\iota(s_{0})}^{\iota(\bar{r})} d \Xi + \sum_{0 < i < j} \mu_{\iota(s^{j})}(\iota(s_{i})) \int_{K} f_{\iota(s_{i})}^{\iota(\bar{r})} d \Xi \\
            & \leq \mu_{\iota(s^{j})}(\iota(s_{0})) (\delta - \varp) + \sum_{0 < i < j} \mu_{\iota(s^{j})}(\iota(s_{i})) \cdot \delta \\
            & = \sum_{ i < j} \mu_{\iota(s^{j})}(\iota(s_{i})) \cdot  \delta - \mu_{\iota(s^{j})}(\iota(s_{0})) \cdot \varp  = \delta - \mu_{\iota(s^{j})}(\iota(s_{0})) \cdot \varp  \leq \delta - \mu(\iota(s_{0})) \cdot \varp.
        \end{split}
    \end{equation*}
    On the other hand, since $\limit_{j \to i^{\ast}} \mu(\iota(s^{j})) = \mu(\iota(r^{\ast}))$ and $\limit_{j \to i^{\ast}} \mu (\iota(r^{\ast}) \sim \iota(s^{j})) = 0,$ there exists $0 < j< \min \{ i^{\ast} + 1, \omega \}$ such that $\mu_{\iota(r^{\ast})}(\iota(r^{\ast}) \sim \iota(s^{j})) < \mu(\iota(s_{0})) \cdot \varp.$ Expressing $\iota(r^{\ast}) = \iota(s^{j}) \vee [\iota(r^{\ast})  \sim \iota(s^{j})]$ and using the previous inequalities, we get: 
    \begin{equation*}
        \begin{split}
            \int_{K} f_{\iota(r^{\ast})}^{\iota(\bar{r})} d \Xi & = \mu_{\iota(r^{\ast})}(\iota(r^{\ast}) \sim \iota(s^{j})) \int_{K} f_{\iota(r^{\ast}) \sim \iota(s^{j})}^{\iota(\bar{r})} d \Xi  +  \mu_{\iota(r^{\ast})}(\iota(s^{j})) \int_{K} f_{\iota(s^{j})}^{\iota(\bar{r})} d \Xi \\
            & < \mu(\iota(s_{0})) \cdot \varp + [\delta - \mu(\iota(s_{0})) \cdot \varp] \\
            & = \delta.
        \end{split}
    \end{equation*}
    Therefore $\int_{K} f_{\iota(r^{\ast})}^{\iota(\bar{r})} d \Xi < \delta,$ so $r^{\ast} \in I.$ Finally, since for all $\ell \in W,$ $  \mu_{\iota(r^{\ast})}(\iota(r_{\ell})) \geq \delta,$ we have that $f_{\iota(r^{\ast})}^{\iota(\bar{r})} \geq \delta,$ hence $\int_{K} f_{\iota(r^{\ast})}^{\iota(\bar{r})}  d \Xi \geq \delta,$ that is, $ r^{\ast} \notin I,$ which is a contradiction. Thus, $I$ is not dense below $r^{\ast}.$ 
\end{PROOF}

Finally, we can prove \autoref{mt}, the main result in this paper.

\begin{theorem}\label{i27}
    Any complete Boolean algebra with a strictly positive probability measure satisfying the $\theta$-density property is $\theta$-$\FAM$-linked, even uniformly $\theta$-$\calY_{**}$-linked where $\calY_{**}\coloneqq \{(\Xi,\bar P)\in \calY_* \colon \Xi \text{ is free}\}$.
\end{theorem}

\begin{PROOF}[\textbf{Proof}]{\autoref{i27}}
    Let $\mu$ a strictly positive probability measure on $\cB$, and fix $S \subseteq \cB^{+}$ witnessing that $\mu$ satisfies the $\theta$-density property. For each  $s \in S$ and $\varp \in (0, 1)_{\bbQ},$ consider  the sets $Q_{s, \varp} \coloneqq \{ b \in \cB \colon \mu_{s}(b) \geq 1 - \varp \}.$ Let us prove that the sequence $\langle Q_{s, \varp} \colon s \in S \conj  \varp \in (0, 1)_{\bbQ} \rangle$ witnesses that $\cB$ is uniformly $\theta$-$\calY_{**}$-linked proving \ref{i4II} of \autoref{i4} and conditions~\ref{i10.1.b} and~\ref{i10.1.c} in \autoref{i10}~\ref{i10.1}. First, for condition~\ref{i4II}, let $(\Xi,\bar P)\in \calY_{**}$, $K\coloneqq K_{\bar P}$, $W\coloneqq W_{\bar P}$, $s \in S,$ $\varp_{0} \in (0, 1)_{\bbQ},$ and $\bar{b} = \langle b_{\ell} \colon \ell \in W \rangle \in {}^{W} \! Q_{s, \varp_{0}}.$ Therefore, for all $\ell \in W,$ $ \mu_{s}(b_{\ell}) \geq 1 - \varp_{0},$ that is, we are under the hypothesis of \autoref{i23},\footnote{Applied to $\bbP \coloneqq \cB^{+}$, $\mu,$ $\iota \colon \bbP \to \cB^{+}$ being the identity function,  $r_\ell \coloneqq b_\ell$ for any $\ell \in W,$ $r^* =s$,  and $\delta \coloneqq 1-\varp_0$\label{i27.f}.} by virtue of which there exists some $b^{\otimes} \leq s$ such that, for all $b \leq b^{\otimes},$ $ \int_{K} f_{b}^{\bar{b}}  d \Xi  \geq 1 - \varp_{0},$ where $f_{b}^{\bar{b}} \colon K \to \bbR$ is defined by $ f_{b}^{\bar{b}}(k) \coloneqq \frac{1}{\vert P_{k} \vert} \sum_{\ell \in P_{k}} \mu_{b}(b_{\ell})$ for any $k \in K$. Thereby, we can define $\limit^{s, \varp_{0}} \colon {}^{W} \! Q_{s, \varp_{0}} \to \cB$ such that $\limit^{s, \varp_{0}}(\bar{b}) \coloneqq b^{\otimes},$ which satisfies that, for any  $b \leq \limit^{s, \varp_{0}}(\bar{b}),$ $ \int_{K} f_{b}^{\bar{b}} d \Xi \geq 1 - \varp_{0}.$
    To show~\ref{i4II}, let $i^{\ast} < \omega$, $s_i\in S$, $\varp_i \in (0,1)_\bbQ$, and $\bar{b}^{i} = \langle b_{\ell}^{i} \colon \ell \in W \rangle \in {}^{W} \! Q_{s_i, \varp_i}$ for each $i < i^{\ast}.$ Let $\varp > 0,$ $F \subseteq K$ finite, and $\langle B_{m} \colon m < m^{\ast} \rangle$ a finite partition of $K$ and $b \in \bbP$ such that $b \leq \limit^{s_{i}, \varp_{i}}(\bar{b}^{i})$ for each $i < i^{\ast}.$ Notice that, by the construction of $\limit^{s_{i}, \varp_{i}},$ for all $i < i^{\ast}$ and $b' \leq b,$ we have that $b' \leq \limit^{s_i, \varp_i}(\bar{b}^{i})$ and therefore, $\int_{K} f_{b'}^{\bar{b}^{i}}  d \Xi \geq 1 - \varp_{i}$ for each $i < i^{\ast}.$ So we are under the hypothesis of \autoref{i19},\footnote{Applied to the same parameters as \Cref{i27.f}.} by virtue of which there are a non-empty  finite set $u \subseteq K \setminus F$ and $b^{\oplus} \leq b$ such that:  

    \begin{enumerate}[label=\rm{(\arabic*)}]
        \item  $\left| \frac{\vert u \cap B_{m} \vert }{\vert u \vert} - \Xi(B_{m}) \right| < \varp,$ for all $m < m^{\ast},$

        \item $\frac{1}{\vert u \vert} \sum_{k \in u} \frac{\vert \{ \ell \in P_{k} \colon b^{\oplus} \leq b_{\ell}^{i} \} \vert}{\vert P_{k} \vert} \geq 1 - \varp_{i} - \varp,$ for all $i < i^{\ast}.$ 
    \end{enumerate}

    
    Second, for condition~\ref{i10.1.b}, let $\varp \in (0, 1)_{\bbQ}$ and $b \in \cB^{+}.$ By the $\theta$-density property, there exists some $s \in S$ such that, $\mu_{s}(b) \geq 1 - \varp,$ that is, $b \in Q_{s, \varp}.$ Thus, not only $\bigcup_{s \in S} Q_{s, \varp_{0}}$ is dense in $\cB$ but it is the whole $\cB$. Finally, condition~\ref{i10.1.c} is a direct consequence of \autoref{s9} (and redundant by \autoref{rm-i10}). 
    %
    %
\end{PROOF}


\begin{remark}\label{i33}
    We expect that \autoref{i19} and variations of \autoref{i23} can be used to find new examples of $\theta$-$\FAM$-linked forcing notions. Concretely, 
    if a forcing notion satisfies statements similar to \autoref{i23} without demanding that $\mu$ is $\sigma$-additive, then \autoref{i19} (or some variation) could be applied as in the proof of \autoref{i27} to show $\theta$-$\FAM$-linkedness for some $\theta$. 
\end{remark}

As a consequence of \autoref{i27}, using \autoref{i3.0}, we get: 

\begin{corollary}\label{i99}
    The measure algebra adding $\theta$-many random reals is $\theta$-$\FAM$-linked (even uniformly $\theta$-$\calY_{**}$-linked) for any infinite cardinal $\theta.$
\end{corollary}

Finally, in the countable case, the Baire $\sigma$-algebra on ${}^{\omega} 2$ coincides with $\calB(\cantor)$. Consequently, $\cB_{\omega}$ is random forcing and therefore, we get a result due to Saharon Shelah (see \cite{Sh00}): 

\begin{corollary}\label{i30}
    Random forcing is $\sigma$-$\FAM$-linked (even uniformly $\sigma$-$\calY_{**}$-linked). 
\end{corollary}

{\small

\begin{thebibliography}{CMUZ25}

\bibitem[BCM25]{BCM2}
J\"{o}rg Brendle, Miguel~A. Cardona, and Diego~A. Mej\'ia.
\newblock Separating cardinal characteristics of the strong measure zero ideal.
\newblock {\em J. Math. Log.}, 2025.
\newblock To appear, \href{https://arxiv.org/abs/2309.01931}{arXiv:2309.01931}.

\bibitem[CMU24]{CMU}
Miguel~A. Cardona, Diego~A. Mejía, and Andrés~F. Uribe-Zapata.
\newblock A general theory of iterated forcing using finitely additive measures.
\newblock Preprint, \href{https://arxiv.org/abs/2406.09978}{arXiv:2406.09978}, 2024.

\bibitem[CMU25]{CMUP}
Miguel~A. Cardona, Diego~A. Mej\'ia, and Andr\'es~F. Uribe-Zapata.
\newblock {F}initely additive measures on {B}oolean algebras.
\newblock Preprint, \href{https://arxiv.org/abs/2503.08910}{arXiv:2503.08910}, 2025.

\bibitem[HS16]{HoroShe}
Haim Horowitz and Saharon Shelah.
\newblock Saccharinity with ccc.
\newblock Preprint, \href{https://arxiv.org/abs/1610.02706}{\texttt{arXiv:1610.02706}}, 2016.

\bibitem[Kel59]{Kelley59}
J.~L. Kelley.
\newblock Measures on {B}oolean algebras.
\newblock {\em Pacific J. Math.}, 9:1165--1177, 1959.

\bibitem[KST19]{KST}
Jakob Kellner, Saharon Shelah, and Anda~R. T\u{a}nasie.
\newblock Another ordering of the ten cardinal characteristics in {C}icho\'{n}'s diagram.
\newblock {\em Comment. Math. Univ. Carolin.}, 60(1):61--95, 2019.

\bibitem[Kun84]{Ku84}
Kenneth Kunen.
\newblock Random and {C}ohen reals.
\newblock In {\em Handbook of set-theoretic topology}, pages 887--911. North-Holland, Amsterdam, 1984.

\bibitem[Kun11]{Kunen}
Kenneth Kunen.
\newblock {\em Set {T}heory}, volume~34 of {\em Studies in Logic (London)}.
\newblock College Publications, London, 2011.

\bibitem[Mej24]{M24Anatomy}
Diego~A. Mejía.
\newblock Anatomy of {$\tilde{\mathbb{E}}$}.
\newblock {\em Ky\={o}to Daigaku S\=urikaiseki Kenky\=usho K\=oky\=uroku}, 2290:43--61, 2024.
\newblock \url{https://www.kurims.kyoto-u.ac.jp/~kyodo/kokyuroku/contents/2290.html}.

\bibitem[MU25]{MU}
Diego~A. Mej\'ia and Andr\'es~F. Uribe-Zapata.
\newblock Probability trees.
\newblock Preprint, \href{https://arxiv.org/abs/2501.07023}{arXiv:2501.07023}, 2025.

\bibitem[Oxt80]{Oxtoby}
John~C. Oxtoby.
\newblock {\em Measure and category}, volume~2 of {\em Graduate Texts in Mathematics}.
\newblock Springer-Verlag, New York-Berlin, second edition, 1980.
\newblock A survey of the analogies between topological and measure spaces.

\bibitem[Ros98]{ross98}
Sheldon~M. Ross.
\newblock {\em A First Course in Probability}.
\newblock Prentice Hall, Upper Saddle River, N.J., fifth edition, 1998.

\bibitem[She00]{Sh00}
Saharon Shelah.
\newblock Covering of the null ideal may have countable cofinality.
\newblock {\em Fund. Math.}, 166(1-2):109--136, 2000.

\bibitem[Uri23]{uribethesis}
Andrés~F. Uribe-Zapata.
\newblock Iterated forcing with finitely additive measures: applications of probability to forcing theory.
\newblock Master's thesis, Universidad Nacional de Colombia, sede Medell\'in, 2023.
\newblock \href{https://www.researchgate.net/publication/369113160_Iterated_forcing_with_finitely_additive_measures_applications_of_probability_to_forcing_theory}{https://sites.google.com/view/andres-uribe-afuz/publications}.

\bibitem[Uri24]{IntNumU24}
Andrés~F. Uribe-Zapata.
\newblock The intersection number for forcing notions.
\newblock {\em Ky\={o}to Daigaku S\=urikaiseki Kenky\=usho K\=oky\=uroku}, 2290:1--17, 2024.
\newblock \url{https://www.kurims.kyoto-u.ac.jp/~kyodo/kokyuroku/contents/2290.html}.

\end{thebibliography}

\bibliographystyle{alpha}
}

\end{document}